\documentclass[12pt, reqno]{amsart}
\usepackage[english]{babel}
\usepackage[utf8]{inputenc}
\usepackage{amsmath, amsthm, amscd, amssymb, fullpage, fancyhdr, upgreek, amssymb, hyperref, enumerate, comment, siunitx, enumitem, graphicx, mathrsfs, mathtools}
\usepackage[dvipsnames]{xcolor}
\usepackage{microtype}
\usepackage[toc]{appendix}
\usepackage[parfill]{parskip}

\newtheorem{theorem}{Theorem}[section]
\newtheorem{proposition}[theorem]{Proposition}

\newtheorem{lemma}[theorem]{Lemma}
\newtheorem{conjecture}[theorem]{Conjecture}

\theoremstyle{remark}

\DeclareMathOperator{\sym}{sym}
\DeclareMathOperator{\USp}{USp}
\DeclareMathOperator{\SpO}{SO}
\DeclareMathOperator{\U}{U}
\DeclareMathOperator{\Prob}{Prob}
\DeclareMathOperator{\ad}{ad}

\DeclareMathOperator{\res}{Res}
\DeclareMathOperator{\vol}{vol}
\DeclareMathOperator{\eff}{eff}

\renewcommand{\mod}{~\mathrm{mod}~}

\usepackage{tikz}
\usepackage{tkz-tab}
\usepackage{tkz-graph}
\usetikzlibrary{shapes.geometric,positioning}

\numberwithin{equation}{section}
\renewcommand{\Im}{\mathrm{Im}}
\renewcommand{\Re}{\mathrm{Re}}
\numberwithin{figure}{section}

\newcommand{\hr}[1]{\href{#1}{\url{#1}}}

\title{A survey of a random matrix model for a family of cusp forms
}

\author{Owen Barrett, Zo\"e X. Batterman, Aditya Jambhale, Steven J. Miller, Akash L. Narayanan, Kishan Sharma, Chris Yao}

\address{Department of Mathematics, University of California, Berkeley, Berkeley, CA 94720}
\email{\textcolor{blue}{\href{mailto:barrett@math.berkeley.edu}{barrett@math.berkeley.edu}}}

\address{Department of Mathematics and Statistics, Pomona College, Claremont, CA 91711}
\email{\textcolor{blue}{\href{mailto:zxba2020@mymail.pomona.edu}{zxba2020@mymail.pomona.edu}}}

\address{Department of Pure Mathematics and Mathematical Statistics, University of Cambridge, Cambridge CB3 0WA, UK}
\email{\textcolor{blue}{\href{mailto:aj644@cam.ac.uk}{aj644@cam.ac.uk}}}

\address{Department of Mathematics and Statistics, Williams College, Williamstown, MA 01267}
\email{\textcolor{blue}{\href{mailto:sjm1@williams.edu}{sjm1@williams.edu}},  \textcolor{blue}{\href{Steven.Miller.MC.96@aya.yale.edu}{Steven.Miller.MC.96@aya.yale.edu}}}

\address{Department of Mathematics, University of Michigan, Ann Arbor, MI 48104}
\email{\textcolor{blue}{\href{mailto:anaray@umich.edu}{anaray@umich.edu}}}

\address{Department of Pure Mathematics and Mathematical Statistics, University of Cambridge, Cambridge CB3 0WA, UK}
\email{\textcolor{blue}{\href{mailto:kds43@cam.ac.uk}{kds43@cam.ac.uk}}}

\address{Department of Mathematics, Yale University, New Haven, CT 06520}
\email{\textcolor{blue}{\href{mailto:chris.yao@yale.edu}{chris.yao@yale.edu}}}

\thanks{This work was completed during the 2023 SMALL REU program at Williams College and is joint with colleagues from earlier summers. It was supported in part by NSF Grants DMS1561945 and DMS1659037, the University of Michigan, and Williams College.}
\date{\today}
\begin{document}
\maketitle

\begin{abstract} 
The Katz-Sarnak philosophy states that statistics of zeros of $L$-function families near the central point as the conductors tend to infinity agree with those of eigenvalues of random matrix ensembles as the matrix size
tends to infinity. 
While numerous results support this conjecture, S. J. Miller observed that for finite conductors, very different behavior can occur for zeros near the central point in elliptic curve families. This led to the excised model of Due\~{n}ez, Huynh, Keating, Miller, and Snaith, whose predictions for quadratic twists of a given elliptic curve are beautifully fit by the data. The key ingredients are relating the
discretization of central values of the $L$-functions to excising matrices based on the value of the characteristic polynomials at 1 and using lower order terms (in statistics such as the one-level density and pair-correlation) to adjust the matrix size. We discuss recent successes by the authors in extending this model to a family of quadratic twists of finite conductor of a given holomorphic cuspidal newform of level an odd prime level. In particular, we predict very little repulsion for forms with weight greater than 2.
\end{abstract}


\section{Introduction}
The Katz-Sarnak philosophy states that statistics of zeros of $L$-function families near the central point as the conductors tend to infinity agree with those of eigenvalues of certain random matrix ensembles as the matrix size
tends to infinity \cite{KS99a, KS99b}. While the general philosophy yields remarkable predictive insights for both local and global statistics, classic matrix ensembles fail to reflect finer statistical properties of $L$-function zeros. 

For instance, in 2006, Miller \cite{Mil06} observed that the elliptic curve $L$-function zero statistics for finite conductors deviated significantly from the scaling limit of the expected model of orthogonal matrices (which was known to be correct in the limit for suitable test functions), though the fit improved as the conductors increased. 
Subsequently, Due\~{n}ez, Huynh, Keating, Miller, and Snaith \cite{DHKMS12} created the \emph{excised orthogonal model} for finite conductors to explain the phenomena Miller observed in the elliptic curve case. The excision of matrices with small values of the characteristic polynomial evaluated at 1 corresponds to the discretization of the central values of normalized elliptic curve $L$-functions at the central point, which by the Kohnen-Zagier formula generalizes and holds for holomorphic cusp forms \cite{Wal80,KZ81,BM07,Mao08}.
For ease of reading, we recall the statement in its most general form—with no restriction on the conductor of the twist—below \cite[Theorem 1.5]{Mao08}.

\begin{theorem}
Let $f$ be a normalized Hecke eigenform of weight $2k$ and odd level $N$, $g$ a Shimura correspondence of $f$, and $L(s,f\otimes\psi_d)$ the $L$-function of $f$ twisted by the quadratic character $\psi_d$ with fundamental discriminant $d$. The formula of Kohnen and Zagier is
\begin{align}\label{eqn:kohnen-zagier-formula}
L(k,f \otimes \psi_d) \ = \ \frac{c(|d|)^2}{|d|^{k-1/2}}\kappa_f, \quad \text{ where } \kappa_f \ = \ \frac{\pi^k}{(k-1)!}\frac{\langle f,f \rangle}{\langle g,g \rangle}\kappa
.\end{align}
\end{theorem}

This observation led Due\~{n}ez, Huynh, Keating, Miller, and Snaith to create a matrix model whose characteristic polynomials evaluated at 1 mimic the values of elliptic curve $L$-functions at the central point. To model the behavior of low-lying zeros using matrices, they find an analogous discretization of the values of the characteristic polynomials at 1 by introducing a \emph{cutoff value}. 
The other key ingredient is a modification of the matrix size of the ensemble. They consider two matrix sizes: one related to the mean density of zeros and the other determined from lower-order terms of the one-level density. The observed data of quadratic twists of fixed elliptic curve $L$-functions showed terrific agreement with these predictions for finite conductors.

\subsection{Overview of the main paper.}

In this survey article, we summarize our results from \cite{BBJMNSY}. 
We work in the general setting of a family of quadratic twists of the $L$-function associated to a cuspidal modular newform $f$ of level an odd prime, arbitrary integral weight, and nebentype $\chi_f$.
We assume the Generalized Riemann Hypothesis for $L$-functions arising from cuspidal newforms. We also assume the Ratios Conjectures, which we state in Section \ref{section:background_notation}, as it allows us to obtain strong estimates on the ratios of logarithmic derivatives. In turn, this yields clearer information on the distribution of the zeros of the $L$-functions in question \cite{CFZ08}. Using results in the literature, we constructed a family of twists of a fixed newform which pivotally depends on the form's nebentype and self-duality. By the modularity theorem, our family for a given form with weight 2 and principal nebentype restricts to the family of elliptic curve $L$-functions of finite conductor considered by Due\~nez, Huynh, Keating, Miller, and Snaith \cite{DHKMS12}.

To create the model, we determined the associated ensemble for the family. 
Assuming the Ratios Conjectures, we derived the lower-order expansion of the scaled one-level density for our family of twists which we then compared to the densities of certain classical compact groups.
In this survey article when we consider self-dual forms with complex multiplication (`self-CM'), we only consider those forms whose $L$-function has $+1$ sign—we made this choice as there are almost no self-CM forms with sign $-1$ on LMFDB. We matched the densities from the random matrix side to the lowest-lying zeros of the family to find the predicted ensemble for each case, listed as follows.
\begin{center}
\begin{tabular}{c|c}
    Case & Group \\
    \hline
    $\chi_f$ principal, even twists & $\SpO(2N)$\\
    $\chi_f$ principal, odd twists & $\SpO(2N+1)$\\
    $\chi_f$ non-principal and $f = \overline{f}$ & $\USp(2N)$\\
    $\chi_f$ non-principal and $f \neq \overline{f}$ & $\U(N)$\\
\end{tabular}
\end{center}

Matching the leading lower-order term allows us to find an effective matrix size for the model for the first three cases, which we state at the end of Section \ref{section:one_level_densities}. 
Since there are no lower-order terms of the one-level density for the unitary ensemble to which we could compare those of our family, we needed to consider another statistic. We turned to pair-correlation to gather arithmetic terms to use for the effective matrix size. Guided by the conjecture in \cite{DHKMS12}, which followed the argument by Bogolomny, Bohigas, Leboeuf, and Monastra in \cite{BBLM06}, we turned to pair-correlation to recover an arithmetic term for the effective matrix size. We again opt to sketch the main ideas in favor of leaving the details to the main paper \cite{BBJMNSY}.

We follow the recipe to determine a cutoff value for the family of even twists of a form with principal nebentype and even weight. The key ingredient to the creation of the cutoff is the Kohnen-Zagier formula in \cite[Theorem 1]{KZ81}, which applies to those forms with level an odd number. Theoretically, we use the heuristic developed in \cite[Section 5]{DHKMS12} and find that the Kohnen-Zagier formula implies a repulsion from the central point for those forms of weight 2 with principal nebentype. For forms with weight 4, we predict very little repulsion from the origin; for even integral weights greater than or equal to 6, we predict no repulsion.
In practice, we discard those characteristic polynomials, that when evaluated at 1, are larger than our proposed cutoff value; a process we call `excision.' 
This excision resolves in the general setting of forms with principal nebentype the repulsion from the central point observed by Miller.
We say an ensemble with the proposed effective matrix size and cutoff value is the \emph{effective excised matrix ensemble}.

\begin{center}
\begin{figure}[htpb]	
	\label{tab:intro-graph}
	\begin{tabular}{c c c}
		Level 5, weight 8 (principal, even) & Level 3, weight 7 (self-CM) \\
  \includegraphics[scale=0.45]{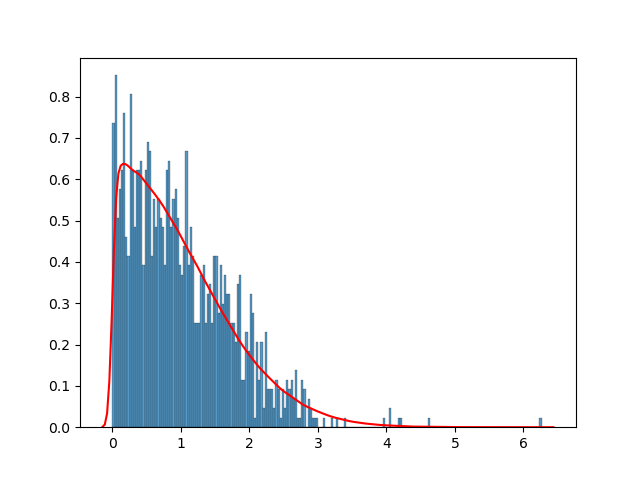} & \includegraphics[scale=0.45]{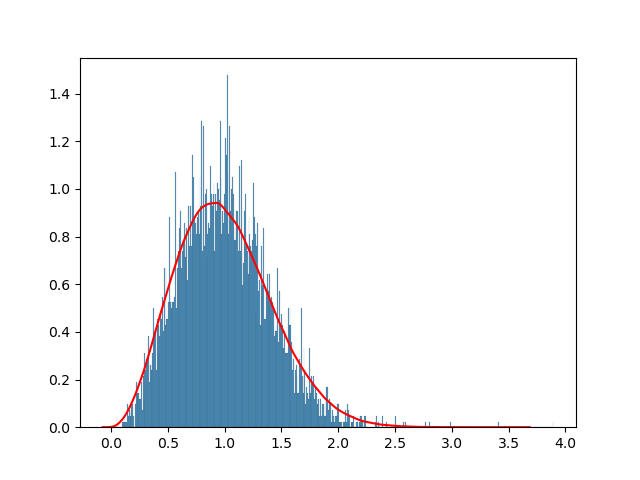}
	\end{tabular}
	\caption{The left histogram shows 1,380 lowest-lying zeros of even twists of \texttt{5.8.a.a} with the red line showing the distribution of the first eigenvalue of 1,000,000 randomly generated special orthogonal matrices with characteristic polynomial evaluated at 1. The right histogram shows the distribution of 5,450 lowest zeros for even twists of \texttt{3.7.b.a} with the red line showing the distribution of first eigenvalues of 1,000,000 randomly generated symplectic matrices with characteristic polynomial evaluated at 1.  The data have been normalized to have mean 1.}
\end{figure}	
\end{center}

In Figure \ref{tab:intro-graph}, we present histograms of the lowest-lying zeros of our family which we compare to the eigenvalues of characteristic polynomials evaluated at 1 of random matrices arising from the random matrix ensembles.

We gather numerical data in Section \ref{sec:numerical-observations} and verify the distribution of lowest-lying zeros of our family matches that of the predicted symmetry groups. We were not able to analytically find the effective matrix size for our family, except for when the form is principal and has weight 2. The difficulty unfortunately boils down to being unable to explicitly determine the Euler product of an arbitrary cuspidal newform $L$-function. 

We verify the statistics of predicted ensembles of a given form aligns with that of the family; we also find a few curious numerics which point us to further theoretical investigations. For instance,
we find that the distribution of zeros of the family associated to certain generic forms (forms with non-principal nebentype and no complex multiplication) recover non-generic behavior (either principal nebentype or self-CM).

\section{Background and notation}\label{section:background_notation}
We collect basic facts about random matrices and $L$-functions to be used in \cite{BBJMNSY}.

\subsection{Matrix ensembles and one-level densities}
\label{sec:matrix-ensembles-one-level-densities}
Denote by $G(N)$ one of $\SpO(2N)$, $\SpO(2N+1)$, $\USp(2N)$, and $\U(N)$. Let $\varphi$ be an even Schwartz function.
Since the mean spacing of the eigenangles of matrices in $G(N)$ depends on $N$, we may scale the eigenangles to have mean spacing one.
We obtain the asymptotic scaled one-level density for large $N$ by scaling the unscaled one-level density found in \cite[Corollary AD.12.5.2]{KS99a} and expanding in powers of either $1/N$ or $1/(2N+1)$. In particular, the scaled one-level density formula is given by 
\begin{align}
R(\varphi,\SpO(2N)) & \ \coloneqq \ \int_0^1 \varphi(\theta)\left(1-\frac{1}{2N} + \frac{\sin\left( (2N-1)\theta\pi /N \right) }{2N \sin(\theta\pi /N)} \right)d\theta\\
& \ = \ \int_0^1 \varphi(\theta)\left(1 + \frac{\sin(2\pi \theta)}{2\pi \theta} - \frac{1 + \cos(2\pi \theta)}{2N}- \frac{\pi \theta \sin(2\pi \theta)}{6N^2} + O(N^{-3})\right)d\theta \nonumber
\end{align}
\begin{align}
R(\varphi,\SpO(2N+1)) & \ \coloneqq \ \int_0^1 \varphi(\theta) \left(1-\frac{1}{2N+1} - \frac{\sin(4\pi\theta N/(2N+1))}{(2N+1)\sin\left(2\pi\theta/(2N+1)\right)}\right)d\theta \label{eq: scaledoneleveldensity_SO2N+1}\\
& \ = \ \int_0^1 \varphi(\theta)\left(1 - \frac{\sin(2\pi\theta)}{2\pi\theta} - \frac{1-\cos(2\pi\theta)}{2N+1} + \frac{2\pi\theta\sin(2\pi\theta)}{3(2N+1)^2} + O(N^{-3})\right)d\theta \nonumber
\end{align}
\begin{align}
R(\varphi,\USp(2N)) & \ \coloneqq \ \int_0^1 \varphi(\theta)\left(1 + \frac{1}{2N} - \frac{\sin\left( (2N+1)\theta \pi /N \right) }{2N \sin(\theta \pi /N)}\right)d\theta \\
& \ = \ \int_0^1\varphi(\theta)\left(1- \frac{\sin (2\pi \theta)}{2\pi \theta} + \frac{1 -\cos(2\pi \theta)}{2N} + \frac{\pi \theta \sin (2\pi \theta)}{6N^2} + O(N^{-3}) \right)d\theta \nonumber\\
R(\varphi,\U(N)) & \ \coloneqq \ \int_0^1 \varphi(\theta)d\theta
.\end{align}

\newpage
\noindent

\subsection[L-functions]{L-functions}\label{sct:Lfunctions} We turn to cuspidal newforms and their associated $L$-functions.

\subsubsection{Cuspidal newforms.}\label{sct:newforms-duals}
We consider the linear space $S_k(N,\chi_f)$ of cusp forms of level $N$, integral weight $k$, and nebentype $\chi_f$ for the Hecke congruence subgroup $\Gamma_0(N)$. We focus on those cuspidal modular forms which are newforms.
In particular, if $f \in S_k^{\text{new}}(M,\chi_f)$, then $f$ is an eigenform and has Fourier expansion $f(z)= \sum_{n = 1}^{\infty} a_f(n) e^{2\pi i nz}$
at the cusp $\infty$ where the $a_f(n)$'s are the Fourier coefficients. Define $\lambda_f(n) \coloneqq a_f(n)n^{-(k-1)/2}$ so that the functional equation of the associated $L$-function relates $s \to 1-s$. 
For $\text{Re}(s) > 1$, the $L$-function $L(s,f)$ associated to $f$ is given by the Dirichlet series $
L(s,f) \ \coloneqq \ \sum_{n\geq1} \lambda_f(n) n^{-s}$
and Euler product
\begin{align}
L(s,f) \ & = \ \prod_p \left( 1- \lambda_f(p)p^{-s} + \chi_f(p) p^{-2s} \right) ^{-1} \nonumber \\
\ & = \ \prod_p \left( 1 - \alpha_f(p)p^{-s} \right) ^{-1} \left( 1 - \beta_f(p) p^{-s} \right) ^{-1},
\end{align}
where the Satake parameters $\alpha_f,\beta_f$ satisfy $\alpha_f(p) + \beta_f(p) = \lambda_f(p)$ and $\alpha_f(p) \beta_f(p) = \chi_f(p)$.
Such an $L$-series admits an analytic continuation to the entire complex plane and has functional equation \cite[Proposition 1.3.6]{Bump97}. By comparing the coefficient of the $p^{-ms}$ term of the Euler product with the Dirichlet series, we obtain the relation
\begin{align} \label{eqn:lambda-coefficient-relation}
\lambda_f(p^{m}) \ = \ \sum_{\ell \geq 0} \alpha_f(p)^{\ell} \beta_f(p)^{m-\ell}
.\end{align}
The form $\overline{f}$ dual to $f$ has Fourier coefficients which satisfy the duality relation $\lambda_f(n) = \chi_f(n) \lambda_{\overline{f}}(n) = \chi_f(n) \overline{\lambda_{f}}(n)$ for $\mathrm{gcd}(n,M)=1$ by the adjointness formula for a cuspidal Hecke form $f$ \cite[Proposition 14.11]{IK04}.
We may relate $L(s,f)$ to $L(1-s,\overline{f})$ by the functional equation
\begin{align}
	L(s,f) \ = \ \epsilon_f \bigg( \frac{\sqrt{M} }{2\pi} \bigg) ^{1-2s} \frac{\Gamma\left( \frac{k+1}{2} - s \right) }{\Gamma\left( \frac{k-1}{2} + s \right) } L(1-s,\overline{f})
,\end{align}
where $\epsilon_f$ is the sign and has absolute value 1.

\subsubsection{Rankin-Selberg convolution}\label{sct:rankin-selberg}
For $f \in S_{k}^{\text{new}}(M_f,\chi_f)$ and $g \in S_{\kappa}^{\text{new}}(M_g,\chi_g)$, the $L$-series of their Rankin-Selberg convolution is defined to be \begin{equation}
L(s,f\otimes g)  \ \coloneqq\ L(2s,\chi_f\chi_g) 
 \sum_{n\geq 1} a_f(n) a_g(n)n^{-s}, \end{equation} provided the least common multiple of $M_f$ and $M_g$ is square-free. It is well-known that $L(s,f\otimes g)$ has analytic continuation and admits an Euler product \cite[Theorem 1.6.2]{Bump97}.
The local factor at an unramified prime (i.e., those primes not dividing the level of $f$ or $g$) is given by
\begin{align}
L_p(s,f\otimes g) & \ = \ 1 - \lambda_f(p)\lambda_g(p)p^{-s} + \left(\chi_f(p)\lambda_g(p)^2 + \chi_g(p)\lambda_f(p)^2 - 2 \chi_f(p) \chi_g(p) \right)p^{-2s} \\
	& \quad \ - \ \chi_f(p)\chi_g(p)\lambda_f(p)\lambda_g(p)p^{-3s} + \chi_f(p)^2\chi_g(p)^2p^{-4s} \nonumber
.\end{align}	
\subsubsection{Quadratic twists}
An integer $d$ is a fundamental discriminant provided that $d$ is either square-free and congruent to 1 modulo 4 or is four times a square-free integer congruent to 2 or 3 modulo 4. Let $L(s,f_d) \coloneqq L(s,f\otimes \psi_d)$ denote the $L$-function obtained by twisting $L(s,f)$ by a quadratic character $\psi_d$ with fundamental discriminant $d$, that is, conductor $|d|$.
For $\mathrm{Re}(s)>1$, the twisted $L$-function $L(s,f_d) \ = \ \sum_{n\geq 1} \lambda_f(n) \psi_d(n)n^{-s}$ has Euler product
\begin{align}\label{eqn:euler-twisted}
L(s,f_d) \ = \ \prod_p \left( 1 - \lambda_f(p)\psi_d(p)p^{-s} + \chi_f(p)\psi_d(p)^2p^{-2s} \right) ^{-1}.
\end{align}
Provided $\mathrm{gcd}(d,M) = 1$, the completed $L$-function satisfies the functional equation\footnote{\cite[Section 14.8]{IK04}.}
\begin{align}
L(s,f_d) \ = \ \epsilon_{f\otimes \psi_d} \bigg( \frac{\sqrt{M} |d|}{2\pi} \bigg) ^{1-2s} \frac{\Gamma\left( \frac{k+1}{2}-s \right) }{\Gamma\left(\frac{k-1}{2} + s \right) } L(1-s,\overline{f}_d)
,\end{align}
with root number $\epsilon_{f\otimes \psi_d} =\epsilon_f\omega_f(d)= \epsilon_f\chi_f(d) \psi_d(M) \tau(\psi_d)^2/d =  \epsilon_f\chi_f(d)\psi_d(-|D|)$ where $\tau(\psi)$ denotes the Gauss sum and the second equality follows from \cite[Corollary 2.1.47]{Coh07} since $\psi_d=\left(\frac{d}{\cdot}\right)$ is a real character.
The approximate functional equation for a twist shifted by $\alpha$ is given in \cite[Section 5.1]{CFZ08} by
\begin{align}\label{eqn:approx-functional-eq-twisted}
L(1 /2+\alpha,f_d) & \ = \ \sum_{m < x} \frac{\lambda_f(m)\psi_d(m)}{m^{1 /2+\alpha}} \\
& \quad \ + \ \omega_f(d) \epsilon_f \bigg( \frac{\sqrt{M} |d|}{2\pi} \bigg) ^{-2\alpha} \frac{\Gamma\left(\frac{k}{2} - \alpha\right)}{\Gamma\left(\frac{k}{2} + \alpha\right)} \sum_{n < y} \frac{\overline{\lambda_f(n)}\psi_d(n)}{n^{1 /2-\alpha}} + \text{remainder}, \nonumber
\end{align}
where $xy=d^2 /2\pi$. We ignore the remainder term \cite{CFKRS05}.
\subsubsection{Complex multiplication and self-duality}\label{sct:CM}
A newform $f$ is said to have complex multiplication (CM) by $\eta$ if there is a non-principal Dirichlet character $\eta$ such that $\eta(p)\lambda_f(p) = \lambda_f(p)$ for all primes $p$ in a set of primes of density 1. We call self-dual newforms with non-principal nebentype constructed by Shimura (see \cite[Section 3 ]{Rib77}) `self-CM' forms. Such forms have complex multiplication by their own nebentype \cite[Remark 2]{Rib77}. We document a property of any self-CM form $f$ constructed by Shimura (c.f. details in \cite{BBJMNSY}):
For each positive fundamental discriminant $d$ prime to $M$, the sign $\epsilon_{f\otimes\psi}$ of the twisted $L$-function is equal to $\epsilon_f$. To prove the equality, it suffices to note that a Dirichlet character whose modulus equals its fundamental discriminant is Kronecker and that twisting in the self-CM case preserves the nebentype \cite[Remark 1]{Rib77}.

We introduce some notation. Let $S_{k}^{\text{new}}(M,\text{principal})$ denote the family of forms with $\chi_f$ principal; let $S_{k}^{\text{new}}(M,\text{self-CM})$ denote the family of self-dual forms with $\chi_f$ non-principal; and let $S_{k}^{\text{new}}(M,\text{generic})$ denote the family of generic forms with $\chi_f$ non-principal and $f \neq \overline{f}$.
\subsubsection{Symmetric and adjoint square L-functions.}\label{sct:sym-adj-sq}
We follow \cite[Section 5.12]{IK04} and define symmetric and adjoint square $L$-functions of newforms as factors of Rankin-Selberg convolutions.
Let $\chi'_f$ denote the primitive character that induces the nebentypus $\chi_f$ of $f$. 
Denote $L(s,\sym^2f)  \ \coloneqq \ L(s,f\otimes f)L(s,\chi'_f)^{-1}$ and $L(s,\ad^2f) \coloneqq L(s,f\otimes \overline{f})\zeta(s)^{-1}$.
The symmetric square $L$-function $L(s,\sym^2f)$ has an Euler product with local factors at unramified primes given by
\begin{align}\label{eqn:local-factor-sym-sq}
	L(s,\sym^2f) = (1-\alpha_f(p)^2p^{-s})^{-1}(1-\alpha_f(p)\beta_f(p)p^{-s})^{-1}(1-\beta_f(p)^2p^{-s})^{-1}
.\end{align}
We record analytic facts about these functions.
First, note that $L(s,\ad^2f)$ is always entire since $L(s,f\otimes \overline{f})$ always has a simple pole at $s=1$ which cancels with the zero of $\zeta(s)$.
For $\mathcal{H}$ the upper half-plane, we have an arithmetically significant value at $s=1$ for $L(s,\ad^2f)$ by the equality 
\begin{align} \label{eqn:residue-at-1-f-bar-f}
	L(1,\ad^2f) \ = \ \frac{(4\pi)^k \langle f, f\rangle}{\Gamma(k)\vol (\Gamma_0(M) \backslash \mathcal{H})} \ = \ \res (L(s,f\otimes\overline{f}),1 )
.\end{align}
When $f \in S_{k}^{\text{new}}(M,\text{principal})$, then $\overline{f} = f$ and $L(s,\ad^2f) = L(s,\sym^2f)$ by $\chi'_f \equiv 1 \mod M$.
When $f$ has non-principal nebentypus, $L(s,\chi_f)$ is entire. The $L$-function $L(s,f\otimes f)$ may not be entire since $f$ may still be self-dual. When $f$ has non-principal nebentypus, then $f$ is self-CM if and only if $\overline{f} = f$ \cite[Section 3, Remark 2]{Rib77}.
If $f$ is self-CM, $L(s,\sym^2f)$ has a pole, and the symmetric and adjoint square $L$-functions do not coincide since $L(s,\chi'_f) \neq \zeta(s)$.
In fact, the case of $f$ self-CM is the only case when $L(s,\sym^2f)$ inherits the pole from $L(s,f\otimes f)$ at $s = 1$ as $L(s,\chi_f')$ is entire.

\subsection{The set of ``good" fundamental discriminants}\label{sct:fundamental-discriminant-family}
We wish to determine the set of ``good" fundamental discriminants for which the ratios recipe (see Section \ref{subsect:ratios_conjectures}) and the Kohnen-Zagier formula may be applied. We discuss in detail the construction of the family in the main paper \cite{BBJMNSY}. In particular, we wish to replace $\epsilon_{f\otimes \psi_d}$ and $\psi_d(-M)$ in the functional equation \eqref{eqn:approx-functional-eq-twisted} to adapt the ratios recipe.

Let $\mathcal{D}$ denote the set of positive fundamental discriminants, and fix a cuspidal newform $f \in S_k^{\text{new}} (M,\chi_f)$. 
Fix an integer $\Delta \in \{+1,-1\}$. The family of fundamental discriminants we wish to study is given by
\begin{align*}
\mathcal{D}_f^{+}(X) \ \coloneqq \ \begin{cases}
\{d \in \mathcal{D} \mid 0<d\leq X,\ \ \psi_d(-M) \epsilon_f \ = \ +1\} & \qquad \text{$\chi_f$ principal, even twists}, \\
\{d \in \mathcal{D} \mid 0<d\leq X,\ \ \psi_d(-M) \epsilon_f \ = \ -1\} & \qquad \text{$\chi_f$ principal, odd twists}, \\
\{d \in \mathcal{D} \mid 0<d\leq X,\ \ \psi_d(-M) \ = \ \Delta\} & \qquad \text{$f$ self-CM}, \\
\{d \in \mathcal{D} \mid 0 < d \leq X\} & \qquad \text{$f$ generic}
.\end{cases}
\end{align*}
Finally, we state estimates on the cardinality $|\mathcal{D}_f^{+}(X)|$ of the family; the detailed calculations will appear in the full paper \cite{BBJMNSY}:
\begin{align}\label{eqn:cardinality-estimates} 
|\mathcal{D}_f^{+}(X)| \ = \ \begin{cases}
3MX(2\pi^2(M+1))^{-1} + O(X^{1 /2}) & \qquad f = \overline{f}, \\
3MX(\pi^2(M^2-1))^{-1} + O(X^{1 /2}) & \qquad f \neq \overline{f}.
\end{cases}
\end{align}

\subsection{The family of quadratic twists}
Denote the family of quadratic twists of a fixed holomorphic cuspidal newform $f$ with fundamental discriminant ranging over $\mathcal{D}^{+}_f(X)$ by $
\mathcal{F}_f^{+}(X)$.
For $d \in \mathcal{D}^{+}_f(X)$, the quadratic character $\psi_d(M) = (d/M)$ assumes the value
\begin{align}
\mathcal{E}_f(M) \ \coloneqq \ \psi_d(M) \ = \ \begin{cases}
-\epsilon_f & \qquad \chi_f \text{ principal, even twists}, \\
\ \ \epsilon_f & \qquad \chi_f \text{ principal, odd twists}, \\
\ \ -\Delta & \qquad \chi_f \text{ non-principal, } f = \overline{f}, \\
(d/M) & \qquad \chi_f \text{ non-principal, } f \neq \overline{f}.
\end{cases}
\end{align}
\subsection{Ratios Conjectures}\label{subsect:ratios_conjectures}
To derive a formula for the one-level density of the zeros near the critical point $1/2$ of $L$-functions for $\mathcal{F}_f^{+}(X)$, we consider the average over the family of ``good" fundamental discriminants of a ratio of shifted $L$-functions:
\begin{align} \label{eqn:ratio-shifted}
R_f(\alpha,\gamma) \ \coloneqq \ \sum_{d \in \mathcal{D}_f^{+}(X)} \frac{L(1/2 + \alpha,f_d)}{L(1/2+ \gamma,f_d)}
.\end{align}
Using~\eqref{eqn:euler-twisted} we have $1/L(s,f_d)=\sum_{n\geq1}\mu_f(n)\psi_d(n)n^{-s}$
where $\mu_f(n)$ is a multiplicative function defined by
\begin{align} \mu_f(n)\ =\ \begin{cases}
  -\lambda_f(p) & \text{if $n=p$,} \\
  \ \ \ \chi_f(p) & \text{if $n=p^2$,} \\
  \ \ \ \ 0 & \text{if $n=p^j$, $j>2$.}
\end{cases} \label{eq:mudef}
\end{align}
We denote the first sum arising from the approximate
functional equation~\eqref{eqn:approx-functional-eq-twisted} by
\begin{align}\label{eq:R1def}R_f^1(\alpha,\gamma)\ \coloneqq\ \sum_{d \in \mathcal{D}_f^{+}(X)}
\sum_{h,m\geq 0}\frac{\lambda_f(m)\mu_f(h)\psi_d(mh)}{m^{1/2+\alpha}h^{1/2+\gamma}}.
\end{align}
Likewise for $R_f^2(\alpha,\gamma)$:
\begin{align}\label{eq:R2def}R_f^2(\alpha,\gamma)\ \coloneqq\ 
\omega_f(d)\epsilon_f\frac{\Gamma\left(\frac k2-\alpha\right)}{\Gamma\left(\frac k2+\alpha\right)}
\sum_{d \in \mathcal{D}_f^{+}(X)}\left(\frac{\sqrt M |d|}{2\pi}\right)^{-2\alpha}
\sum_{h,m\geq 0}\frac{\overline{\lambda_f(m)}\mu_f(h)\psi_d(mh)}
{m^{1/2-\alpha}h^{1/2+\gamma}}.
\end{align}
Our formulation of the Ratios Conjectures is stated as follows:\footnote{For other examples of the Ratios Conjectures, see \cite{HMM11}.}
\begin{conjecture}\label{conj:ratiosconjecture}
For the conditions $-1/4 < \mathrm{Re}(\alpha) < 1/4$, $1/\log x \ll \mathrm{Re}(\gamma) < 1/4$ and $\mathrm{Im}(\alpha),\Im(\gamma) \ll X^{1-\varepsilon}$, the average over the family of a ratio of shifted $L$-functions is
\begin{align}\label{eqn:R_f(alpha,gamma)-error}
R_f(\alpha,\gamma) & \ = \ \sum_{d\in \mathcal{D}_f (X)} \left[ Y_f A_f(\alpha,\gamma) + \eta_f \bigg( \frac{\sqrt{M} |d|}{2\pi} \bigg)^{-2\alpha}  \frac{\Gamma\left(k/2 - \alpha\right)}{\Gamma\left(k/2 + \alpha\right)} \widetilde{Y}_f \widetilde{A}_f (-\alpha,\gamma) \right] \\
& \qquad \qquad \qquad \qquad \qquad \qquad \qquad \qquad \qquad \qquad \qquad + \ O(X^{1/2+\varepsilon}), \nonumber
\end{align}
where 
\begin{align}
Y_f(\alpha,\gamma) & \ = \ \frac{L(1 + 2\gamma,\chi'_f)L(1+2\alpha,\sym^2f)}{L(1+\alpha+\gamma,\chi'_f)L(1+\alpha+\gamma,\sym^2f)}, \label{eqn:Y_f(alpha,gamma)} \\ 
\widetilde{Y}_f(-\alpha,\gamma) & \ = \ \frac{L(1+2\gamma,\chi'_f) L(1-2\alpha,\sym^2\overline{f})}{\zeta(1-\alpha+\gamma)L(1-\alpha+\gamma,\ad^2f)}, \label{eqn:tilde-Y_f(alpha,gamma)}
\end{align}
\begin{align} \label{eqn:A_f(alpha,gamma)}
A_f(\alpha,\gamma) & \ = \ Y_f(\alpha,\gamma)^{-1} V_{\mid}(\alpha,\gamma) V_{\nmid}(\alpha,\gamma)
,\end{align}
\begin{align}\label{eqn:V_mid(alpha,gamma)}
V_{\mid}(\alpha,\gamma) \ = \ \sum_{m=0}^{\infty} \left( \frac{\lambda_f(M^{m}) \mathcal{E}_f^{m}(M)}{M^{m(1 /2+\alpha)}} - \frac{\lambda_f(M)}{M^{1 /2 + \gamma}} \frac{\lambda_f(M^{m}) \mathcal{E}_f^{m+1}(M)}{M^{m(1 /2+\alpha)}} \right)
,\end{align}
\begin{align} \label{eqn:V_nmid(alpha,gamma)}
V_{\nmid}(\alpha,\gamma) \ = \ \prod_{p\nmid M} \left( 1 + \frac{p}{p+1} \left( \sum_{m=1}^{\infty} \frac{\lambda_f(p^{2m})}{p^{m(1+2\alpha)}}- \frac{\lambda_f(p)}{p^{1+\alpha+\gamma}}\sum_{m=0}^{\infty} \frac{\lambda_f(p^{2m+1})}{p^{m(1+2\alpha)}} + \frac{\chi_f(p)}{p^{1+2\alpha}} \sum_{m=0}^{\infty} \frac{\lambda_f(p^{2m})}{p^{m(1+2\alpha)}} \right)  \right),
\end{align}
\begin{align} \label{eqn:tilde-A_f(alpha,gamma)}
\widetilde{A}_f(-\alpha,\gamma) & = \widetilde{Y}_f(-\alpha,\gamma)^{-1} \widetilde{V}_{\mid}(-\alpha,\gamma) \widetilde{V}_{\nmid}(-\alpha,\gamma),
\end{align}
\begin{align}\label{eqn:widetilde_V_mid(alpha,gamma)}
\widetilde{V}_{\mid}(-\alpha,\gamma) \ = \ \sum_{m=0}^{\infty} \left( \frac{\overline{\lambda_f(p^{m})} \mathcal{E}_f^{m}(p)}{p^{m(1 /2-\alpha)}} - \frac{\lambda_f(p)}{p^{1 /2 + \gamma}} \frac{\overline{\lambda_f(p^{m})} \mathcal{E}_f^{m+1}(p)}{p^{m(1 /2-\alpha)}} \right) ,
\end{align}
\begin{align}\label{eqn:widetilde_V_nmid(alpha,gamma)}
\widetilde{V}_{\nmid}(-\alpha,\gamma) = \prod_{p\nmid M} \left( 1 + \frac{p}{p+1} \left( \sum_{m=1}^{\infty} \frac{\overline{\lambda_f(p^{2m})}}{p^{m(1-2\alpha)}}- \frac{\lambda_f(p)}{p^{1-\alpha+\gamma}}\sum_{m=0}^{\infty} \frac{\overline{\lambda_f(p^{2m+1})}}{p^{m(1-2\alpha)}} + \frac{\chi_f(p)}{p^{1-2\alpha}} \sum_{m=0}^{\infty} \frac{\overline{\lambda_f(p^{2m})}}{p^{m(1-2\alpha)}} \right)  \right),
\end{align}
and the expectation of $\omega_f(d)\epsilon_f$ over $d$ is
\begin{align}\label{eq:rootexpect}
\eta_f \ \coloneqq \ \langle\omega_f(d)\epsilon_f\rangle \ = \ \begin{cases}
+1 & \qquad \text{$\chi_f$ principal, even twists}, \\
-1 & \qquad \text{$\chi_f$ principal, odd twists}, \\
+1 & \qquad \chi_f \text{ non-principal, $f = \overline{f}$},\\ 
\langle\omega_f(d)\epsilon_f\rangle & \qquad \text{$\chi_f$ non-principal, $f \neq \overline{f}$}.
\end{cases}
\end{align}
\end{conjecture}
The authors in \cite{HKS09} follow the recipe of \cite{CFKRS05}, \cite{CFZ08} and the calculations in \cite{CS07} to derive a formula for \eqref{eqn:R_f(alpha,gamma)-error}.
Our version of the Ratios Conjectures is the same as Conjecture 2.1 in \cite{HKS09} with appropriate substitutions for our family where we use the corresponding approximate functional equation \eqref{eqn:approx-functional-eq-twisted}. 
In particular, the definition \eqref{eq:mudef} implies we need only consider $h=0,1$
for~\eqref{eqn:V_mid(alpha,gamma)} and $h=0,1,2$ for~\eqref{eqn:V_nmid(alpha,gamma)}.
For our family, we factor out the divergent part of $R^1_f(\alpha,\gamma)$ using the Dirichlet $L$-function $L(s,\chi_f')$ and also the symmetric square $L$-function associated with $L(s,f)$.
The bounds at $1/4$ allow us to control the convergence of
Euler products of the type~\eqref{eqn:V_nmid(alpha,gamma)}. We similarly obtain $R^2_f(\alpha,\gamma)$.
We conclude with our desired formula. We leave the detailed computations to the main paper \cite{BBJMNSY} in favor of sketching the main ideas.

\section{One-level density}\label{section:one_level_densities}

\subsection{Averaging the logarithmic derivative}
To calculate the one-level density, we need the average of the logarithmic derivative of $L$-functions defined by
\begin{align}
\sum_{d \in \mathcal{D}_f^{+}(X)} \frac{L'}{L}(1 /2 + r, f_d) \ = \ \frac{\partial}{\partial \alpha} \bigg|_{\alpha=\gamma=r} R_f (\alpha,\gamma).
\end{align}
The equality follows from differentiating Equation \eqref{eqn:ratio-shifted} and setting $\alpha = \gamma = r$.
\begin{proposition}\label{prop:avg-logarithmic-derivative} Assume the Ratios Conjectures and that $1 /\log X \ll \mathrm{Re}(r) < 1 /4$ and $\mathrm{Im}(r) \ll X^{1- \varepsilon}$.
Then the average of the logarithmic derivative over the family $\mathcal{F}_f^{+}(X)$ is
\begin{align}
\sum_{d\in \mathcal{D}_f^{+}(X)} \frac{L'}{L} (1 /2 + r, f_d) & \ = \ \sum_{d\in \mathcal{D}_f^{+}(X)} \bigg( - \frac{L'}{L}(1+2r,\chi'_f) + \frac{L'}{L}(1+2r,\sym^2f) + A^{1}_f(r,r) \\		
& \quad \ - \ \eta_f \bigg( \frac{\sqrt{M} |d|}{2\pi} \bigg)^{-2r} \frac{\Gamma\left( \frac{k}{2} - r\right)}{\Gamma\left( \frac{k}{2} + r\right)} \frac{L(1+2r,\chi'_f)L(1-2r,\sym^2\overline{f})}{L(1,\ad^2f)} \nonumber \\
& \qquad \ \times \ \widetilde{A}_f (-r,r) \bigg) + O(X^{1/2 + \varepsilon}) \nonumber
\end{align}
where $\widetilde{A}_f$, $A^{1}_f(r,r)$ are given by \eqref{eqn:tilde-A_f(alpha,gamma)} and $
A^{1}_f(r,r)= \frac{\partial}{\partial \alpha} \bigg|_{\alpha=\gamma=r} A_f(\alpha,\gamma)$, respectively.
\end{proposition}

\begin{proof}
	The proof is the same as that of Theorem 2.2 in \cite{HKS09}, mutatis mutandis.
\end{proof}

\subsection{Unscaled one-level density}
With a formula for the average of logarithmic derivatives, we turn to finding the lower-order terms of the scaled one-level density functions for our family.
For some Schwartz function $\varphi$ and $\gamma_d$ ordinates of the zeros on the critical line, the scaled one-level density function is defined by $D_1(\varphi,f) \coloneqq \sum_{d \in \mathcal{D}_f^{+}(X)} \sum_{\gamma_d} \varphi(\gamma_d)$.
We apply the argument principle to rewrite $D_1(\varphi,f)$ as
\begin{align}
\sum_{d \in \mathcal{D}_f^{+}(X)} \frac{1}{2\pi i} \left( \int_{(c)} - \int_{(1-c)} \right) \frac{L'}{L}(s,f_d) \varphi(i(1 /2-s ))\,ds
,\end{align}
where $1 /2 + 1 /\log X < c < 3 /4$ is fixed and $(c)$ denotes the path from $c-i\infty$ to $c+i\infty$.
Turning to the integral on the $c$-line, we have
\begin{align}\label{eqn:unscaled-integrand}
\frac{1}{2\pi} \int_{\mathbb{R}} \varphi(t-i(c-1 /2)) \sum_{d \in \mathcal{D}_f^{+}(X)} \frac{L'}{L} (c + it,f_d) \, dt,\end{align}
and the sum over $d$ can be replaced by Proposition \ref{prop:avg-logarithmic-derivative}.
The bounds on the size $t$ coming from the Ratios Conjectures do not pose a problem; see \cite{HKS09,CS07}.
\begin{lemma}\label{lem:unscaled-regular-integrand}
The integrand in~\eqref{eqn:unscaled-integrand} is regular at $t = 0$.
\end{lemma}

\begin{proof} Since we assumed the Generalized Riemann Hypothesis for Dirichlet $L$-functions, $L(s,\chi_f)$ does not vanish on the line $\mathrm{Re}(s)=1$.
In particular, the poles of $L(s,\chi_f')$ and $L(s,\sym^2f)$ at $s=1$ may obstruct regularity of the integrand in~\eqref{eqn:unscaled-integrand} at $t = 0$.
The $L$-function $L(s,\chi_f)$ has a pole at $s = 1$ if and only if $\chi_f$ is principal, and $L(s,\sym^2f)$ has a pole only when $f$ is self-CM.

For the case $f \in S_{k}^{\text{new}}(M,\text{principal})$, $\widetilde{A}_f(-r,r) = A_f(-r,r)$ which implies $A_f$ is analytic.
As $r \to 0$, $\zeta(1+2r) = (2r)^{-1} + O(1)$ and $-(\zeta' /\zeta)(1+2r) = (2r)^{-1} + O(1)$.
After substituting the expansions into the expression in Proposition \ref{prop:avg-logarithmic-derivative}, the integrand in~\eqref{eqn:unscaled-integrand} is regular at $t=0$ if and only if $\widetilde{A}_f(0,0)=1$, which happens since $A_f(r,r)=1$.
		
For the case $f \in S_k^{\text{new}}(M,\text{self-CM})$, the symmetric square $L$-function of $f$ inherits a pole at $s=1$ from $L(s,f\otimes f)$, and $A_f$ is analytic. Using the relation $\lambda(p^{2m+1})\lambda(p) = \lambda(p^{2m+2}) + \lambda(p^{2m})$ for $p \nmid M$ and the multiplicativity of $\lambda(p)$ for $p = M$ (see Section \ref{sct:newforms-duals}), we get that $\widetilde{A}_f(-r,r) = A_f(-r,r)=1$. We use the assumption $\eta_f=+1$ to obtain
\begin{align}
\lim_{r \to 0} \bigg[ & \frac{L'}{L} (1+2r,\sym^2f) - \left( \frac{\sqrt{M} |d|}{2\pi} \right) ^{-2r} \frac{\Gamma\left(\frac{k}{2}  - r \right)}{\Gamma\left(\frac{k}{2} + r \right)} \frac{L(1+2r,\chi'_f) L(1 -2r, \sym^2\overline{f})}{L(1,\ad^2f)} \widetilde{A}_f (-r,r) \bigg] \nonumber\\
& = \ \lim_{r \to 0} \left( -(2r)^{-1} + O(1) \right) - \frac{L(1+2r,\chi'_f) L(1-2r,\chi'_f)^{-1} L(1-2r,f\otimes f)}{\res(L(s,f\otimes f), 1)} \nonumber\\   
& = \ \lim_{r\to 0}\left( - (2r)^{-1} + O(1) \right) - \left( -(2r)^{-1} + O(1) \right) \  = \ O(1), 
\end{align}
which completes the proof.
\end{proof}	

\begin{proposition}\label{prop:unscaled-density}
Assume the Generalized Riemann Hypothesis and the Ratios Conjectures. The unscaled one-level density for the zeros of the family $\mathcal{F}_f^{+}(X)$ is
\begin{align}\label{eqn:unscaled-density}
D_1(\varphi,f) \ = \ \frac{1}{2\pi} & \int_{\mathbb{R}} \varphi(t) \sum_{d \in \mathcal{D}_f^{+}(X)} \bigg(\log \bigg( \frac{\sqrt{M} |d|}{2\pi} \bigg)^2 + \frac{\Gamma'}{\Gamma}\bigg( \frac{k}{2} + it \bigg) + \frac{\Gamma'}{\Gamma} \bigg( \frac{k}{2} - it \bigg) \\
& \qquad + \ \frac{L'}{L} \bigg( \frac{1}{2} + it, f_d \bigg) + \frac{L'}{L} \bigg( \frac{1}{2}+ it, \overline{f}_d \bigg) \bigg) \, dt + O(X^{1/2 + \varepsilon}) \nonumber
.\end{align}
\end{proposition}	

\begin{proof}	

By Lemma~\ref{lem:unscaled-regular-integrand} we may move the path of integration to $c = 1 /2$, and the rest of the proof is the same as that of Theorem 2.3 in \cite{HKS09}, mutatis mutandis.
\end{proof}

\subsection{Scaled one-level density} 
We scale the lowest-lying zeros to have mean spacing one, and we calculate the one-level density for the scaled zeros. The goal is to recover the next-to-leading order term from equation \eqref{eqn:unscaled-density}.
Motivated by \cite[Theorem 5.8]{IK04}, we rescale the variable $t$ by $\tau = \frac{tR}{\pi}$ where we equate the mean densities of eigenvalues with the mean densities of zeros by setting\footnote{\cite{DHKMS12}}
\begin{align}
R \ \coloneqq \ 
\begin{dcases} \log \bigg( \frac{\sqrt{M} X}{2\pi} \bigg) - \frac{1}{2} & \qquad \chi_f \text{ principal}, \text{odd twists}, \\
\log \bigg( \frac{\sqrt{M} X}{2\pi} \bigg) & \qquad \chi_f \text{ principal, even twists or $f$ self-CM},\\
2 \log \bigg( \frac{\sqrt{M} X}{2\pi} \bigg) & \qquad f \neq \overline{f},
\end{dcases}
\end{align}
and define the even test function $g(\tau) \coloneqq \varphi(t)$. We define the scaled one-level density function to be $S_1(g,f) \coloneqq D_1(\varphi,f)/|\mathcal{D}_f^{+}(X)|$.
Summation by parts and the cardinality estimate \eqref{eqn:cardinality-estimates} yield the following approximations:
\begin{align}
\sum_{d\in \mathcal{D}_f^{+}(X)} \log \bigg( \frac{\sqrt{M} |d|}{2\pi} \bigg) \ & = \ |\mathcal{D}_f^{+}(X)| \bigg( \log\bigg( \frac{\sqrt{M} X}{2\pi }\bigg)-1 \bigg)   + O(X^{1 /2}), \\
\sum_{d\in \mathcal{D}_f^{+}(X)}\bigg( \frac{\sqrt{M} |d|}{2\pi} \bigg) ^{-2i\pi\tau /R} \ & = \ |\mathcal{D}_f^{+}(X)| \bigg( 1 -  \frac{2i \pi\tau}{R} \bigg)^{-1} e^{-2i\pi\tau} + O(X^{1 /2})
.\end{align}
We now wish to obtain a series expansion of the scaled one-level density in terms of $R$.
\begin{proposition}
Assume the Ratios Conjecture. 
Then the scaled one-level density for the zeros of the family $\mathcal{F}_f^{+}(X)$ is given by
\begin{align}
S_1(g,f) \ = \ \int_{\mathbb{R}} g(\tau) \left( 1 + Q(\tau) + O(R^{-3}) \right) d\tau,
\end{align}
where the lower order terms not in the error term are
\begin{align}
Q(\tau)\ =\  \begin{dcases}
\frac{\sin(2\pi \tau)}{2\pi\tau} - a_1 \frac{1 + \cos(2\pi \tau)}{R} - a_2 \frac{\pi \tau \sin(2\pi \tau)}{R^{2}}& \quad \chi_f \ \text{\emph{principal}}, \text{\emph{ even twists}}, \\
- \frac{\sin(2\pi\tau)}{2i\pi\tau} - a_3 \frac{1-\cos(2\pi\tau)}{2R+1} + a_4 \frac{2\pi\tau\sin(2\pi\tau)}{(2R+1)^2} & \quad \chi_f \ \text{\emph{principal}}, \text{\emph{ odd twists}}, \\
- \frac{\sin(2\pi \tau)}{2\pi\tau}  + b_1 \frac{1 - \cos(2\pi \tau)}{R} + b_2 \frac{\pi \tau \sin(2\pi \tau)}{R^2} & \quad \chi_f \ \text{\emph{non-principal}},\ f = \overline{f}, \\ 
\frac{c_1 + c_2 \cos(2\pi \tau)}{R} + d_1 \frac{\pi \tau \sin(2\pi \tau)}{R^2}  & \quad \chi_f \ \text{\emph{non-principal}},\ f \neq \overline{f}
\end{dcases}
\end{align}
with coefficients 
\begin{align}
a_1 \ & = \ 1 - \psi(k/2) - A^{1}(0,0) + \gamma - \frac{L'}{L}(1,\sym^2f) \label{eqn:coeff-a_1}, \\
a_2 \ & = \ -2\psi(k/2) - 2\psi(k/2)\gamma + 2\gamma - 2\gamma_1 + (2\psi(k/2)-2-2\gamma-B'(0)) \frac{L'}{L}(1,\sym^2f) \label{eqn:coeff-a_2} \\
& \quad + \ (\gamma + 1 - \psi(k/2)) B'(0) + \frac{1}{4}B''(0) + 2\frac{L''}{L}(1,\sym^2f) , \nonumber\\
a_3 \ & = \ 2 - 2\psi(k/2) + 2\gamma_1 - 2 \frac{L'}{L}(1,\sym^2f) - 2 A^1(0,0), \label{eqn:coeff-a_3}\\
a_4 \ & = \ 4\psi(k /2) - 4i\pi\tau \gamma + 4\psi(k /2) \gamma + 4\gamma_1 +(2\psi(k /2) -2 - 2\gamma) B'(0) \label{eqn:coeff-a_4} \\
& + ( 4 + 4\gamma + 2B'(0) - 4\psi(k /2))  \frac{L'}{L}(1,\sym^2f) - \frac{B''(0)}{2} - \frac{L''}{L}(1,\sym^2f), \nonumber\\
b_1 \ & = \ 1-\psi(k/2) - \xi_0 \frac{L(1,\chi'_f)}{L(1,\ad^2f)} - A^{1}_f(0,0) + \frac{L'}{L}(1,\chi'_f) \label{eqn:coeff-b_1} ,\\
b_2 \ & = \ -2\psi(k/2) +B'(0)- \psi(k/2) B'(0)+ \frac{B''(0)}{4} + 2 \frac{L''}{L}(1,\chi'_f) \label{eqn:coeff-b_2}\\
&+\frac{L'}{L}(1,\chi'_f)\Big(-2\xi_0 +B'(0) +2-2\psi(\textstyle{\frac{k}{2}})\Big) + \frac{L(1,\chi'_f)}{L(1,\ad^2f)}\Big(2\psi(\textstyle{\frac{k}{2}}) \xi_0 -2\xi_0 +2\xi_1 -\xi_0 B'(0)\Big) ,\nonumber
\end{align}
\begin{align}
c_1 \ & = \ \psi(k/2) + \frac{1}{2} ( (A^{1}_f + A_{\overline{f}}^{1})(0,0)  - \frac{L'}{L}(1,\chi'_f) - \frac{L'}{L}(1,\chi'_{\overline{f}}) + \frac{L'}{L}(1,\sym^2f) + \frac{L'}{L}(1,\sym^2\overline{f}) ) \label{eqn:coeff-c_1} ,\\
c_2 \ & = \ -\frac{1}{2} ( \eta_f \widetilde{A}_f(0,0)L(1,\chi'_f) \frac{L (1,\sym^2\overline{f})}{L(1,\ad^2f)} + \eta_{\overline{f}} \widetilde{A}_{\overline{f}}(0,0)L(1,\chi'_{\overline{f}}) \frac{L(1,\sym^2f)}{L(1,\ad^2\overline{f})} )  \label{eqn:coeff-c_2},\\
d_1 \ & = \ \eta_f \frac{L(1,\sym^2\overline{f})}{L(1,\ad^2f)}\bigg( -\frac{1}{2}\widetilde{B}'_f(0)L(1,\chi'_f) + \psi(k /2) \widetilde{A}_f(0,0)L(1,\chi'_f) - \widetilde{A}_f(0,0)L(1,\chi'_f) \bigg) \label{eqn:coeff-d_1} \nonumber \\
& \quad + \ \eta_{\overline{f}} \frac{L(1,\sym^2f)}{L(1,\ad^2\overline{f})} \bigg( -\frac{1}{2}\widetilde{B}'_{\overline{f}}(0)L(1,\chi'_{\overline{f}}) + \psi(k /2) \widetilde{A}_{\overline{f}}(0,0)L(1,\chi'_{\overline{f}}) - \widetilde{A}_{\overline{f}}(0,0)L(1,\chi'_{\overline{f}}) \bigg) \nonumber\\
& \quad  + \ \eta_f \frac{L'(1,\sym^2\overline{f})}{L(1,\ad^2f)} \widetilde{A}_f(0,0) L(1,\chi'_f) + \eta_f \frac{L'(1,\sym^2f)}{L(1,\ad^2\overline{f})} \widetilde{A}_{\overline{f}}(0,0) L(1,\chi'_{\overline{f}}),
\end{align}
 where $\psi \coloneqq \Gamma'/{\Gamma}$ is the digamma function $B_f(s) = A_f(-s,s)$, $-B'(0)/2 = A^{1}(0,0)$, and $B_f^{(n)}(s) = \frac{d^n}{d r^n} \Big |_{r=s} A_f(-r,r)$.
\end{proposition}
\begin{proof}
The details of the proof are left to the main paper \cite{BBJMNSY}. We note important observations which allowed us to find scaled one-level density for our family.
For $f \in S_{k}^{\text{new}}(M,\text{principal})$, we get the simplifications $\widetilde{A}_f = A_f$ and $L(1+r,\chi'_f) = \zeta(1+r)$. The proof is the same as that of the expansion (3.19) in \cite{HKS09} with appropriate substitutions.

For $f \in S_k^{\text{new}}(M,\text{self-CM})$ with $+1$ sign, $\eta_f = +1$ and $L(s,\sym^2f)$ has a simple pole at $s=1$ with residue $L(1,\ad^2 f)/L(1,\chi'_f)$ by its factorization.
We find the Laurent expansion at $s=0$ of $L(1+s,\sym^2f)$ and the logarithmic derivative $(L'/L)(1+s,\sym^2f)$ evaluated at $1+s$.
Since $\chi_f$ is non-principal, then $L(s,\chi_f')$ has no pole, and we may write it and its logarithmic derivative as Taylor expansions. 
We substitute the expansions in the rescaled version of \eqref{eqn:unscaled-density} and clear out odd terms to obtain our desired result.

For $f$ generic, $L(s,\chi'_f)$ and $L(s,\sym^2f)$ are entire (see Section \ref{sct:sym-adj-sq}), and the values of $L(s,\chi'_f)$ and its derivatives at $s=1$ are not well known except in particular cases.
Removing the odd terms in the integrand above as they integrate to zero, we obtain our desired result.
\end{proof}

From the computations above, the one-level density of the newforms converges to that of the one-level scaled density of eigenvalues near 1 in certain compact groups. 
We choose the effective matrix size as follows:
\begin{align}
N_{\eff} \ = \ 
\begin{dcases}
\log \bigg( \frac{\sqrt{M} X}{2\pi} \bigg)/2a_1 & \qquad \text{$\chi_f$ principal, even twists},\\
\left(\log \bigg( \frac{\sqrt{M} X}{2\pi} \bigg) - \frac{1}{2} \right) /a_3-1/2 & \qquad \text{$\chi_f$ principal, odd twists},\\
2 \log \bigg( \frac{\sqrt{M} X}{2\pi} \bigg)/2b_1 & \qquad \text{$f$ self-CM},
\end{dcases}
\end{align}
to match the leading lower-order term.  
Recall from Section \ref{sec:matrix-ensembles-one-level-densities} that there are no lower-order terms for the scaled one-level density expansion of unitary matrices. Hence, we cannot match the leading lower-order term in the generic case to choose an effective matrix size for that case. The authors in \cite{DHKMS12} follow the argument made by Bogomolny, Bohigas, Leboeuf, and Monastra in \cite{BBLM06} to conjecture that the improvement made by using matrices of size $N_{\eff}$ also holds for all $n$-point correlation or density functions. The conjecture motivates the use of the pair-correlation statistic to find the effective matrix size in the generic case.

\section{Pair-correlation}\label{section:pair_correlation}

The series expansion of the scaled one-level density for the unitary group has no lower-order terms, and so we cannot refine the fit of the random matrix model for the generic case.
In particular, there is no way to incorporate into the model the arithmetic information coming from the leading lower-order term of the scaled one-level density in the generic case. 

When comparing the series expansion for the scaled one-level density in the generic case with the series expansion for the unitary group coming from random matrix theory, we see that the lower-order terms cannot be matched. 
Since we cannot extract any arithmetic in the one-level density of the generic case, we turn to pair-correlation. Pair-correlation is insensitive to any finite set of zeros and is computed for one $L$-function only. We obtain a series expansion for the pair-correlation statistics of a single $L$-function in large $T$ to obtain the lower-order terms of arithmetic origin necessary for the effective matrix size.

Let $\gamma,\gamma'$ denote the imaginary coordinates of non-trivial zeros of $L(s,f_d)$, and suppose $\varphi(s)$ is a holomorphic function throughout the strip $|\mathrm{Im}(s)| < 2$, real-valued on the real line, even, and satisfies the bound $\varphi(x) \ll 1 /(1+x^2)$ as $x \to \infty$. Since we assume the Ratios Conjectures, we do not require that the Fourier transform $\widehat{\varphi}$ is compactly supported or decays rapidly, as is the case in \cite{RS96}. Throughout the argument, we assume the GRH. We would like to evaluate the pair-correlation statistic $P(f_d; \varphi) \coloneqq \sum_{0 < \gamma,\gamma' < T} \varphi(\gamma - \gamma')$.
We start with the formula for the average of the logarithmic derivative of shifted $L$-functions:
\begin{align}\label{eqn:avg-log-derivative-shifted-pair}
	\int_0^{T} \frac{L'}{L}(s+\alpha,f_d) \frac{L'}{L}(1-s+\beta,\overline{f}_d)\,dt,
\end{align}
where $s = 1/2 + it$.
We apply Conjecture 5.1 in \cite{CFZ08} to
\begin{align}
	\mathcal{T}_{f_d}(\alpha,\beta,\gamma,\delta) \ \coloneqq \ \int_{0}^{T} \frac{L(s+\alpha,f_d) L(1-s+\beta,\overline{f}_d)}{L(s+\gamma,f_d)L(1-s+\delta,\overline{f}_d)} \,dt,
\end{align}
since $\mathcal{F}_f$ for $\chi_f$ non-principal and $f \neq \overline{f}$ has unitary symmetry.
Hence, we substitute $K=L=1$ and the group $\Xi_{1,1} = \{(1),(12)\} $ which consists of the identity permutation and the transposition $(12)$, and we identify $\alpha_1 = \alpha$, $\alpha_2 = -\beta$, $\gamma_1=\gamma$, and $\delta_1 = \delta$ into the Ratios Conjectures.
We now state the Ratios Conjecture for our family.
\begin{conjecture}\label{conj:ratios-lemma}
	For $-1 /4 < \mathrm{Re}(\alpha), \mathrm{Re}(\beta) <1 /4$, $1/\log(T) \ll \mathrm{Re}(\delta) < 1/4$, and $\mathrm{Im}(\alpha), \mathrm{Im}(\beta) \ll_{\varepsilon} T^{1-\varepsilon}$ for all $\varepsilon > 0$, we have
	\begin{align}
		\mathcal{T}_{f_d}(\alpha,\beta,\gamma,\delta) \ & = \ \int_0^{T} \bigg( Y_U(\alpha,\beta,\gamma,\delta) A_L(\alpha,\beta,\gamma,\delta) \nonumber \\
		& \quad + \ \bigg( \frac{\sqrt{M} |d| t}{2\pi} \bigg)^{-2(\alpha+\beta)}  Y_U(-\beta,-\alpha,\gamma,\delta) A_L(-\beta,-\alpha,\gamma,\delta) \bigg) dt + O(T^{1 /2+\varepsilon}),
	\end{align}
	where
	\begin{align}
		Y_U(\alpha,\beta,\gamma,\delta) \ = \ \frac{L(1+\alpha+\beta,f_d\otimes \overline{f}_d)L(1+\gamma+\delta,f_d\otimes \overline{f}_d)}{L(1+\alpha+\delta,f_d\otimes \overline{f}_d)L(1+\beta+\gamma,f_d\otimes \overline{f}_d)}
,	\end{align}
	and
	\begin{align}
		A_{L}(\alpha,\beta,\gamma,\delta) \ & = \ \prod_p \frac{(1-p^{-(1+\alpha+\beta)})(1-p^{-(1+\gamma+\delta)})}{(1-p^{-(1+\alpha+\delta)})(1-p^{-(1+\beta+\gamma)})} \nonumber \\
		  & \quad \times \ \sum_{m+h=n+k}\frac{\mu\lambda}{p^{(1 /2+\alpha)m+(1 /2+\beta)n+(1 /2+\gamma)h+(1 /2+\delta)k}}
	.\end{align}
Here, $\mu = \mu_d$ is the coefficient on $n^{-s}$ of the reciprocal series $L(s, f_d)^{-1}$ for $\mathrm{Re}(s) > 1$; explicitly,
\begin{equation}
    \mu(n) \ = \ 
    \begin{cases}
        \lambda(n) & n = p ,\\
        \chi(n) & n = p^2 ,\\
        0 & n = p^j \text{ for } j > 2.
    \end{cases}
\end{equation}
\end{conjecture}
We also remark that
\begin{equation}
    A_L(\alpha,\beta,\gamma,\delta) \ = \ Y_U(\alpha,\beta,\gamma,\delta)^{-1} W_{\mid}(\alpha,\beta,\gamma,\delta) W_{\nmid}(\alpha,\beta,\gamma,\delta),
\end{equation}
where
\begin{align}
  &W_\mid(\alpha,\beta,\gamma,\delta) \ = \ \prod_{p\mid M}
    \begin{multlined}[t]\Bigg[
      \sum_{m=0}^{\infty}\frac{\left|\lambda(p^m)\right|^2}{p^{(1+\alpha+\beta)m}}
      -\sum_{m=0}^{\infty}\frac{\lambda(p^{m+1})\overline\lambda(p^m)\overline\lambda(p)}
      {p^{(1+\alpha+\beta)m+(1+\alpha+\delta)}} \\
      -\sum_{m=0}^{\infty}\frac{\lambda(p^m)\overline\lambda(p^{m+1})\lambda(p)}
      {p^{(1+\alpha+\beta)m+(1+\beta+\gamma)}}
      +\sum_{m=0}^{\infty}\frac{\left|\lambda(p^m)\right|^2\left|\lambda(p)\right|^2}
      {p^{(1+\alpha+\beta)m+(1+\gamma+\delta)}}
      \Bigg],
    \end{multlined} \\
  \intertext{and}
  &W_\nmid(\alpha,\beta,\gamma,\delta) \ = \ \prod_{p\nmid N}
    \begin{aligned}[t]\Bigg[
      &\sum_{m=0}^{\infty}\frac{\left|\lambda(p^m)\right|^2}{p^{(1+\alpha+\beta)m}}
      -\sum_{m=0}^{\infty}\frac{\lambda(p^{m+1})\overline\lambda(p^m)\overline\lambda(p)}
      {p^{(1+\alpha+\beta)m+(1+\alpha+\delta)}}
      -\sum_{m=0}^{\infty}\frac{\lambda(p^m)\overline\lambda(p^{m+1})\lambda(p)}
      {p^{(1+\alpha+\beta)m+(1+\beta+\gamma)}} \\
      &+\sum_{m=0}^{\infty}\frac{\overline\lambda(p^m)\lambda(p^{m+2})\overline\chi(p)}
      {p^{(1+\alpha+\beta)m+(2+2\alpha+2\delta)}}
      +\sum_{m=0}^{\infty}\frac{\lambda(p^m)\overline\lambda(p^{m+2})\chi(p)}
      {p^{(1+\alpha+\beta)m+(2+2\beta+2\gamma)}} \\
      &+\sum_{m=0}^{\infty}\frac{\left|\lambda(p^m)\right|^2\left|\lambda(p)\right|^2}
      {p^{(1+\alpha+\beta)m+(1+\gamma+\delta)}}
      -\sum_{m=0}^{\infty}\frac{\lambda(p^{m+1})\overline\lambda(p^m)\lambda(p)\overline\chi(p)}
      {p^{(1+\alpha+\beta)m+(2+\alpha+\gamma+2\delta)}} \\
      &-\sum_{m=0}^{\infty}\frac{\lambda(p^m)\overline\lambda(p^{m+1})\overline\lambda(p)\chi(p)}
      {p^{(1+\alpha+\beta)m+(2+\beta+2\gamma+\delta)}}
      +\sum_{m=0}^{\infty}\frac{\left|\lambda(p^m)\right|^2\left|\chi(p)\right|^2}
      {p^{(1+\alpha+\beta)m+(2+2\gamma+2\delta)}}
      \Bigg].
    \end{aligned}
\end{align}
We used the definition of $\mu$ to state that, if $p \mid N = M |d|^2$, we are free to discard all terms except $h,k \in \{0,1\}$. Otherwise, we may discard all terms except $h,k \in \{0,1,2\}$. See the main paper for details \cite{BBJMNSY}.

\subsection{Averaging the logarithmic derivative} 
We obtain the formula for \eqref{eqn:avg-log-derivative-shifted-pair} by differentiating the result of Conjecture \ref{conj:ratios-lemma} which allows us to compute $P(f_d; \varphi)$ using contour integration. By expanding in series the formula for $P(f_d; \varphi)$ in large $T$, we obtain the desired lower-order terms. Again, we relegate the details to the main paper \cite{BBJMNSY}.
\begin{proposition}\label{prop:int-unramified}
Assume Conjecture \ref{conj:ratios-lemma}, and let $\alpha$, $\beta$, $\gamma$, and $\delta$ be as above. Then
\begin{align}
    & \int_0^{T} \frac{L'}{L}(s+\alpha,f_d) \frac{L'}{L}(1-s+\beta,\overline{f}_d) \, dt \ = \ \int_0^{T} \bigg( \bigg( \frac{L'_{\star}}{L_{\star}} \bigg)' (1+\alpha+\beta,f_d\otimes \overline{f}_d) \nonumber \\
	& \qquad + \ \frac{1}{c_{f_d}^2}\bigg( \frac{\sqrt{M} |d|t}{2\pi} \bigg) ^{-2(\alpha+\beta)} L(1+\alpha+\beta,f_d\otimes \overline{f}_d)L(1-\alpha-\beta,f_d\otimes \overline{f}_d) \nonumber \\
	& \qquad \times \ A_L(-\beta,-\alpha,\alpha,\beta) + \mathscr{C}(1+\alpha+\beta) \bigg)\,dt + O(T^{1 /2+\varepsilon})
.\end{align}
Here, $N \coloneqq M |d|^2$, $c_{f_d} = \emph{Res}(L(s, f_d \otimes \overline{f}_d), 1)$, $\mathscr{C}(1+\alpha+\beta)$ is given in \eqref{script_C}, and $L_{\star}(s, f_d \otimes \overline{f}_d)$ is the unramified part of the Rankin-Selberg convolution as defined in \eqref{rankin_selberg}.
\end{proposition}
\begin{proof}
For notation, put $\Lambda \coloneqq Y_U(\alpha,\beta,\gamma,\delta)
A_L(\alpha,\beta,\gamma,\delta)$ and 
\begin{align}
\Omega \coloneqq \bigg( \frac{\sqrt{M}|d|t}{2\pi} \bigg)^{-2(\alpha+\beta)} Y_U(-\beta,-\alpha,\gamma,\delta) A_L(-\beta,-\alpha,\gamma,\delta).\nonumber
\end{align}
We turn to the derivative of $\Omega$. By applying the usual product rule and setting $\gamma \mapsto \alpha$ and $\beta \mapsto \delta$ by taking the limits appropriately, we are left with
\begin{align}
\frac{\partial^2\Omega}{\partial\beta\,\partial\alpha} \bigg|_{(\gamma,\delta) \ = \ (\alpha,\beta)} \ = \ &\bigg( \frac{\sqrt{M}|d| t}{2\pi} \bigg)^{-2(\alpha+\beta)} \bigg(\frac{L'}{L^2} \bigg)^2(1, f_d\otimes \overline{f}_d) \\
&\times \ L(1-\alpha-\beta, f_d\otimes \overline{f}_d)L(1+\alpha+\beta, f_d\otimes \overline{f}_d)A_L(-\beta,\alpha,\alpha,\beta) \nonumber.
\end{align}
We now wish to evaluate $(L'/L^2)^2(1, f_d \otimes \overline{f}_d)$. Since $L(1, f_d \otimes \overline{f}_d)$ has a simple pole at 1 and $L'(1,f_d\otimes \overline{f}_d)$ has a double pole at 1, the logarithmic derivative $(L'/L)(s, f_d \otimes \overline{f}_d)$ has a simple pole at 1. 
To compute its residue, take a small circular contour $\mathcal{C}$ oriented counterclockwise with center at $s = 1$. By the residue theorem and the argument principle, we have
\begin{equation}
    \text{Res}\bigg( \frac{L'}{L}(s, f_d \otimes \overline{f}_d), 1 \bigg) \ = \ \frac{1}{2\pi i} \oint_{\mathcal{C}} \frac{L'}{L}(s, f_d \otimes \overline{f}_d)\,ds \ = \ \#\{\text{zeros of } L\} - \#\{\text{poles of } L\}.
\end{equation}
The contour may be chosen sufficiently small so that $L(1, f_d \otimes \overline{f}_d)$ has no zeros within the contour which gives the residue $\res( (L'/L)(s, f_d \otimes \overline{f}_d), 1 ) =  -1$. 
The residue of $L(s, f_d \otimes \overline{f}_d)$ at $s = 1$ is $c_{f_d}$ (see Equation \eqref{eqn:residue-at-1-f-bar-f}), and so $(L'/L^2)(s, f_d \otimes \overline{f}_d)$ is entire with value $-1/c_{f_d}$ at $s=1$.
We get
\begin{align}
    \frac{\partial^2\Omega}{\partial\beta\,\partial\alpha} \bigg|_{(\gamma,\delta) \ = \ (\alpha,\beta)} \ &= \ \frac{1}{c_{f_d}^2} \bigg( \frac{\sqrt{M}|d| t}{2\pi} \bigg)^{-2(\alpha+\beta)} L(1-\alpha-\beta, f_d \otimes \overline{f}_d) \\
    & \quad \times \ L(1+\alpha+\beta, f_d \otimes \overline{f}_d)A_L(-\beta,\alpha,\alpha,\beta) \nonumber.
\end{align}
We turn to the derivative of $\Lambda$ evaluated at $(\gamma,\delta) = (\alpha,\beta)$. Then
\begin{align}\label{lambda_derivative_step_1}
    & \frac{\partial^2\Lambda}{\partial\beta\,\partial\alpha} \bigg|_{(\gamma,\delta) \ = \ (\alpha,\beta)} =
    \frac{L''(1+\alpha+\beta, f_d \otimes \overline{f}_d)}{L(1+\alpha+\beta, f_d \otimes \overline{f}_d)}
    -\bigg(\frac{L'(1+\alpha+\beta, f_d \otimes \overline{f}_d)}{L(1+\alpha+\beta, f_d \otimes \overline{f}_d)}\bigg)^2\\
    & \qquad \qquad - \left( \frac{L'}{L} \right)'(1+\alpha+\beta,f_d \otimes \overline{f}_d) + \frac{\partial^2}{\partial\beta\,\partial\alpha}\bigg|_{(\gamma,\delta) = (\alpha,\beta)} (W_{\mid}(\alpha,\beta,\gamma,\delta) + W_{\nmid}(\alpha,\beta,\gamma,\delta)), \nonumber
\end{align}
since $A_L(\alpha,\beta,\alpha,\beta) = Y_U(\alpha,\beta,\alpha,\beta) = 1$.
We focus our attention to the derivative of $W_{\mid}(\alpha,\beta,\gamma,\delta)$ and fix a prime $p$ dividing $M$, which by assumption, must be $M$.
Noting the Fourier coefficients $\lambda_{f_d}(p)$ are completely multiplicative at primes dividing $M$, we use the formula for the sum of a geometric series to obtain
\begin{equation}\label{eqn:second-partial-W-mid-log-sq}
    \frac{\partial^2}{\partial\beta\,\partial\alpha}\bigg|_{(\gamma,\delta) \ =  \ (\alpha,\beta)} W_{\mid}(\alpha,\beta,\gamma,\delta) \ = \ \frac{\log(M)^2}{ |\lambda_{f_d}(M)|^{-2} \cdot M^{1+\alpha+\beta} - 1}.
\end{equation}
We turn to the derivative of $W_{\nmid}(\alpha,\beta,\gamma,\delta)$. The computations are considerably more involved as there is no guarantee the Fourier coefficients $\lambda_{f_d}(p)$ are multiplicative. The strategy involves carefully re-indexing the sums appearing in $W_{\nmid}(\alpha,\beta,\gamma,\delta)$ and repeatedly applying the relation \eqref{eqn:lambda-coefficient-relation}.
We then expand $\lambda_{f_d}(p^{m+2}) - \chi(p) \lambda_{f_d}(p^m)$ in terms of the Satake parameters attached to $f_d$ (see Section \ref{sct:newforms-duals}).
Letting $L_{\star}$ denote the unramified part of $L$, given by the Euler product
\begin{equation}\label{rankin_selberg}
L_{\star}(s) \ = \ \prod_{p\nmid N} L_p(s) \ = \ L(s) \prod_{p \mid N} L_p(s)^{-1}
,\end{equation}
in the half-plane of convergence and its analytic continuation elsewhere, we obtain
\begin{align*}
& \frac{\partial^2}{\partial\beta\,\partial\alpha}\bigg|_{(\gamma,\delta)= (\alpha,\beta)} W_{\nmid}(\alpha,\beta,\gamma,\delta) = \bigg( \frac{L_{\star}'}{L_{\star}} \bigg)'(1+\alpha+\beta, f_d \otimes \overline{f}_d) - \sum_{p \nmid N} \bigg[ \bigg( \frac{\alpha_{f_d}(p) \overline{\alpha}_{f_d}(p) \log(p)}{p^{1+\alpha+\beta} - \alpha_{f_d}(p) \overline{\alpha}_{f_d}(p)} \bigg)^2 \\
& + \bigg( \frac{\alpha_{f_d}(p) \overline{\beta}_{f_d}(p) \log(p)}{p^{1+\alpha+\beta} - \alpha_{f_d}(p) \overline{\beta}_{f_d}(p)} \bigg)^2 + \ \bigg( \frac{\overline{\alpha}_{f_d}(p) \beta_{f_d}(p) \log(p)}{p^{1+\alpha+\beta} - \overline{\alpha}_{f_d}(p) \beta_{f_d}(p)} \bigg)^2 + \bigg( \frac{\beta_{f_d}(p) \overline{\beta}_{f_d}(p) \log(p)}{p^{1+\alpha+\beta} - \beta_{f_d}(p) \overline{\beta}_{f_d}(p)} \bigg)^2\bigg].\nonumber
\end{align*}
Let $\mathscr{B}(1+\alpha+\beta)$ denote the sum over $p\nmid N$ on the RHS.
Note that the Ramanujan-Petersson conjecture is known for holomorphic cusp forms of even integral weight. In particular, we get that $|\alpha_{f_d}| = |\beta_{f_d}| = 1$ for primes $p$ not dividing the level $N$; see \cite{IK04} and \cite{Sar05}. Therefore, $\mathscr{B}(v)$ exists and is analytic in a neighborhood of $v = 1$. Finally, put
\begin{equation}\label{script_C}
    \mathscr{C}(1 + \alpha + \beta) \ \coloneqq \ -\mathscr{B}(1 + \alpha + \beta) + \frac{\log(M)^2}{ |\lambda_{f_d}(M)|^{-2} \cdot M^{1+\alpha+\beta} - 1}.
\end{equation}
Putting everything together yields the result.
\end{proof}

\subsection{Contour integration for pair-correlation}
The formula for the average of the logarithmic derivative for shifted $L$-functions allows us to obtain a formula for $P(f_d; \varphi)$, which will later be used to perform a series expansion to obtain lower-order terms of arithmetic origin with which to calibrate our effective matrix size.
Set
\begin{align}
I_r \ \coloneqq \ &\int_{-T}^{T} \varphi(r) \bigg( 2\log^2\bigg( \frac{\sqrt{M} |d| t}{2\pi} \bigg) + \bigg( \frac{L'_{\star}}{L_{\star}} \bigg)' (1 + ir, f_d \otimes \overline{f}_d) \\
&+ \ \frac{1}{c_{f_d}^2} \bigg( \frac{\sqrt{M} |d| t}{2\pi} \bigg)^{-2ir} L(1+ir, f_d \otimes \overline{f}_d)L(1-ir, f_d \otimes \overline{f}_d) \mathscr{A}(ir) + \mathscr{C}(1+ir) \bigg)\,dr \nonumber.
\end{align}
\begin{proposition}\label{contour_theorem}
Assuming the Ratios Conjectures and with $\varphi$ as above, we have
\begin{equation}
    P(f_d; \varphi) = \frac{1}{2\pi^2} \int_{0}^{T} \bigg[ 2\pi\varphi(0) \log\bigg( \frac{\sqrt{M}|d| t}{2\pi} \bigg) + I_r \bigg]\,dt + O(T^{1/2 + \varepsilon}).
\end{equation}
Here, $I_r$ should be regarded as a principal-value integral near $r=0$.
\end{proposition}
The proof is that of \cite[Section 4]{CS07}, mutatis mutandis, and mainly relies on contour integration and some asymptotic analysis; the formula for \eqref{eqn:avg-log-derivative-shifted-pair} will help with the integrals. 

\subsection{Series expansion and effective matrix size}\label{subsect:series_pair_correlation}
With a formula for $P(f_d; \varphi)$, we obtain a series development for large $T$ and use it to obtain the effective matrix size.
We scale the pair-correlation by substituting $y \coloneqq rR/\pi$ and $R \coloneqq \log \bigg( \frac{\sqrt{M} |d| T}{2\pi e} \bigg)$.
Define the \emph{rescaled test function} $g$ by $g(y) = g(rR/\pi) \coloneqq \varphi(r)$.
By Proposition \ref{contour_theorem} and changing variables within the integral via $r \mapsto rR/\pi= y$, we observe
\begin{equation}\label{changeofvariable}
    \sum_{0 < \gamma,\gamma' < T} g\bigg((\gamma-\gamma')\frac{R}{\pi}\bigg) \ = \ \frac{1}{2\pi^2} \int_{0}^{T} \bigg[ 2\pi g(0) \log\bigg( \frac{\sqrt{M}|d| t}{2\pi} \bigg) + I_y\bigg]\,dt + O(T^{1/2 + \varepsilon}),
\end{equation}
where
\begin{align}
    I_y \ \coloneqq \ &\frac{\pi}{R} \int_{-T(R/\pi)}^{T(R/\pi)} g(y) \bigg( 2\log^2\bigg( \frac{\sqrt{M} |d| t}{2\pi} \bigg) + \bigg( \frac{L'_{\star}}{L_{\star}} \bigg)' \bigg(1 + \frac{i\pi y}{R}, f_d \otimes \overline{f}_d\bigg) \\
    &+ \ \frac{1}{c_{f_d}^2} \bigg( \frac{\sqrt{M} |d| t}{2\pi} \bigg)^{-2i\pi y/R} L\bigg(1+\frac{i\pi y}{R}, f_d \otimes \overline{f}_d\bigg)L\bigg(1-\frac{i\pi y}{R}, f_d \otimes \overline{f}_d\bigg) \nonumber \\
    & \quad \times \ \mathscr{A}\bigg(\frac{i\pi y}{R}\bigg) + \mathscr{C}\bigg(1+\frac{i\pi y}{R}\bigg) \bigg)\,dy \nonumber.
\end{align}
Recall that the bound for test function $\varphi(x) \ll 1/(1+x^2)$ for real $x$ implies a bound for the rescaled test function $g(y)\ll 1/(1+y^2)$ for real $y$. The bound allows us to approximate the pair-correlation statistic for large $T$ by
\begin{align}\label{beforeexpansion}
& \frac{T}{\pi}\log \bigg( \frac{\sqrt{M} |d| T}{2\pi e} \bigg) \bigg[ g(0) + \int_{\mathbb{R}} g(y) \bigg( \frac{1+R^2}{R^2} + \frac{1}{2R^2} \bigg(\frac{L_{\star}'}{L_{\star}}\bigg)' \bigg(1+\frac{i\pi y}{R}, f_d \otimes \overline{f}_d\bigg) \nonumber\\
&\qquad + \ \frac{1}{2R^2c_{f_d}^2} \frac{e^{-2\pi i y(1 + 1/R)}}{1 - 2\pi i y /R} L\bigg(1+\frac{i\pi y}{R}, f_d \otimes \overline{f}_d\bigg)L\bigg(1-\frac{i\pi y}{R}, f_d \otimes \overline{f}_d\bigg) \mathscr{A}\bigg(\frac{i\pi y}{R}\bigg) \nonumber\\
&\qquad+ \ \frac{1}{2R^2} \mathscr{C}\bigg(1+\frac{i\pi y}{R}\bigg)\bigg)\,dy \bigg
] + O(T^{\varepsilon+1/2}).
\end{align}
In the main paper \cite{BBJMNSY}, we performed a series expansion in $1/R$ for $R$ large. Most of the computations are routine, but critical to the argument is the fact that the Ramanujan-Petersson conjecture is known for holomorphic cusp forms. This allows us to expand $\mathscr{B}$ in a series about 1.

After defining the function $h(y) \coloneqq 2 g(y)$ (and thus $h$ also satisfies $h(y) \ll 1/(1+y^2)$), we get \begin{align}\label{afterexpansion}
P(f_d;\varphi) \ = \ &\frac{T}{2\pi}\log \bigg( \frac{\sqrt{M} |d| T}{2\pi e} \bigg) \bigg[ h(0) + \int_{\mathbb{R}} h(y) \bigg( 1 - \bigg(\frac{\sin\pi y}{\pi y}\bigg)^2 \\
&+ \frac{e_1 - e_2 \sin^2 \pi y}{R^2} - \frac{e_3 \pi y \sin 2\pi y}{R^3} + O(R^{-4}) \bigg)\,dy \bigg] + O(T^{\varepsilon+1/2}), \nonumber
\end{align}
where
\begin{align}
    e_1 \ &\coloneqq \ \frac{1}{2} \cdot \frac{\log(M)^2}{ |\lambda_{f_d}(M)|^{-2} \cdot M - 1} , \qquad e_3 \ \coloneqq \ \frac{16 + \mathscr{A}'''(0)}{12}, \\
    e_2 \ &\coloneqq \ -2 + \gamma^2 + 2\gamma_1 - \frac{\mathscr{A}''(0)}{2} - \bigg(\frac{L'}{L}\bigg)'(1, \text{ad}^2\,f_d),
\end{align}
and $\gamma$ is the Euler-Mascheroni constant and $\gamma_1$ is the first Stieltjes constant. We remark that this agrees with and explicitly extends Montgomery's pair-correlation conjecture in \cite{Mon73}; refer to the main paper \cite{BBJMNSY} for details.

Now that we obtained a series expansion for the pair-correlation statistic, we can find the effective matrix size in the generic case. In particular, we obtain the effective matrix size for the entire family $\mathcal{F}_f^{+}(X)$ simply by averaging over all choices $d$ of quadratic twists. 
This fact simplifies the argument as we do not need to average over an infinite family of $L$-functions.
Following \cite{Con05}, the scaled pair-correlation $Q_{\U(N)}(x)$ for $U(N)$ is
\begin{equation}
    Q_{\U(N)}(x) \ = \ 1 - \bigg( \frac{\sin \pi x}{\pi x} \bigg)^2 - \frac{\sin^2 \pi x}{3N^2} + O(N^{-4}).
\end{equation}
Due to the presence of the $e_1$ term in our pair-correlation expansion, we are unable to directly match coefficients to obtain the effective matrix size. To deal with this, we minimize the function $(e_1 - e_2 \sin^2 \pi y)/R^2 + \sin^2 \pi y/3N^2$
with respect to the $L^2$ norm and obtain the effective matrix $N =R/\sqrt{3e_2 - 4e_1}$ for the $L$-function attached to a fixed twisted form $f_d$. For our family, we take the average value of $e_1$ and $e_2$ as $d$ varies over the family and obtain the effective matrix size
\begin{equation}
    N_{\text{eff}} \ = \ \frac{R}{\sqrt{3\langle e_2 \rangle - 4 \langle e_1 \rangle}}
,\end{equation}
for the entire family $\mathcal{F}_f^{+}(X)$ provided that $3\langle e_2 \rangle - 4 \langle e_1 \rangle > 0$.
\section{Cutoff Value}\label{appendix:cutoff_asymptotics}
We optimize the cutoff value of the excision threshold by following \cite{CKRS05} and \cite{CKRS06} as modified by \cite{DHKMS12}.
For our family of even twists of a given form with principal nebentype and even integral weight, we may apply Kohnen-Zagier's formula \cite{Mao08} to get
\begin{align}
	L(1 /2,f_d) \ < \ \frac{|c(d)|^2\kappa_f}{|d|^{k-1/2}} \ \implies \ L(1 /2,f_d) \ = \ 0
,\end{align}
where the coefficients $c(d)$ are the Fourier coefficients of a half-integral weight modular form that are obtained via the generalized Shimura correspondence.
In particular, the arithmetic of the coefficients $c(d)$ is not well-understood.
Though a notion of Shimura correspondence exists for forms with non-principal nebentype, there is no analogous formula of Waldspurger type which means we cannot predict the coefficient of the main term of the frequency of vanishing. For this reason, it is not known whether the values $L(1 /2,f_d)$ are discretized for forms with non-principal nebentype. In the main paper \cite{BBJMNSY}, we provide numerical evidence that the lowest-lying zeros of family of the generic form \texttt{13.2.ea} show (minimal) repulsion at the origin. 

As in \cite{DHKMS12,KeSn00}, for $\Re(s)>-1/2$, the moment generating function $M_O(2N,s)$ of the values $|\Lambda_A(1)|$ as $A$ varies in the random matrix ensemble $\SpO(\text{even})$ can be evaluated as
\begin{align}
M_O(2N,s) = \int_{\SpO(2N)}\Lambda_A(1)^s\,dA=2^{2Ns}\prod_{j=1}^N \frac{\Gamma(N+j-1)\Gamma(s+j-1/2)}{\Gamma(j-1/2)\Gamma(s+j+N-1)}
.\end{align}
For a constant $c > 0$, the probability density for values of the characteristic polynomials $\Lambda_A(1)$ with $A \in \SpO(2N)$ is given by
\begin{align}
P_O(2N,x) = \frac{1}{2\pi ix}\int_{c-i\infty}^{c+i\infty}M_O(2N,s)x^{-s}ds
.\end{align}
We may predict the frequency of vanishing of $\mathcal{F}_f^{+}(X)$ for $f$ with principal nebentypus, even twists by calculating the probability that a random variable $Y_d$ with probability density 
\begin{align}
P_{f}(d,x) \coloneqq \frac{1}{2\pi i x} \int_{c-i\infty}^{c+i\infty} a_f(s) M_O(s,\log d) x^{-s}ds \sim a_{f}(-1/2)P_O(s,\log d)
\end{align}
assumes a value less than $|c(d)|^2\kappa_f|d|^{k-1/2}$ and then sum asymptotically over the family.
This method was pioneered in \cite{CKRS05} and \cite{CKRS06}.
However, to obtain an asymptotic for the frequency of vanishing with the correct leading coefficient, we need to have an idea of the statistical behavior of the coefficient $c(d)$ evaluated over $\mathcal{D}_f^{+}(X)$.

For this reason the authors in \cite{DHKMS12} determine this value numerically; that is, they introduce a notion of an `effective' cutoff depending on a parameter that is determined numerically, but which does not depend on $d$.
We define $\delta_f \kappa_f |d|^{k-1/2}$ for our `effective' cutoff where $\delta_f$ is a numerical input. In particular, our model requires only the parameter $\delta_f$ be numerically determined as in \cite[Section 5]{DHKMS12}—no further inputs are necessary. We write
\begin{align}
\Prob\bigg( 0 \leq Y_d \ \leq \ \frac{\delta_f \kappa_f}{|d|^{k-1/2}} \bigg) \ & \sim \ \int_0^{\delta_f\kappa_f |d|^{1/2-k}} a_f(-1 /2) \ h(\log d) \ x^{-1 /2} \ dx \\	
 & = \ 2 a_f(-1 /2) \ h(\log d)\ \frac{\sqrt{\delta_f \kappa_f}}{|d|^{k/2 - 1/4}}\nonumber
,\end{align}
where $h(N) = 2^{-7 /8}G(1 /2)\pi^{-1 /4}(2N)^{3 /8}$ is the asymptotic for the moment generating function of the symmetry group for $\mathcal{F}_f^{+}(X)$ at the pole and $G$ is the Barnes $G$-function \cite{Bar00}. Following \cite{CKRS05} and \cite{DHKMS12}, we predict that
\begin{align}\label{eqn:conj-starred-sum}
	\# \bigg\{ L_f(s,\psi_d) \in \mathcal{F}_f^{+}(X) : d \text{ prime}, \ L_f(1 /2,\psi_d) = 0 \bigg\}  \ = \ {{\sum_{\substack{d\in \mathcal{D}_f^{+}(X) \\ d \text{ prime}}}}}\Prob\bigg( 0\leq Y_d \leq \frac{\delta_f\kappa_f}{|d|^{k-1/2}} \bigg).
\end{align}
There are asymptotically $X /4\log X$ prime fundamental discriminants $d \in \mathcal{D}_f^{+}(X)$.
The convergence of the sum \eqref{eqn:conj-starred-sum} is determined by $k$. Namely, if $k< 3$, then the sum diverges, and we have the asymptotic
\begin{align}
\eqref{eqn:conj-starred-sum} \ & \sim \ 
\frac{1}{4\log X}2a_f(-1 /2) \sqrt{\delta_f\kappa_f} h(\log X) \frac{4}{5-2k} |X|^{(5-2k) /4} .
\end{align}
For $k=1$, we may find $\sqrt{\delta_f}$ by numerically determining the left hand side of \eqref{eqn:conj-starred-sum}. These numerics have already been run in \cite{DHKMS12}. When $k=2$, there is a finite number of even twists that vanish which indicates there is little to no repulsion of lowest-lying zeros at the origin. When $k \geq 3$, then the sum \eqref{eqn:conj-starred-sum} converges, indicating that at most a finite number of lowest-lying zeros vanish \cite{CKRS05}; hence, no repulsion at the origin is expected. In Section \ref{sec:numerical-observations}, we verify these phenomena numerically for the forms: \texttt{11.2.a.a}, \texttt{7.4.a.a}, \texttt{5.8.a.a}.
\section{Numerical Observations}
\label{sec:numerical-observations}
We reiterate that we did not compute the effective matrix size for our family due to the given cuspidal newform $L$-function's Euler product not being accessible in PARI/GP or SageMath.
Recall the nearest neighbour spacing statistic is the probability density for distances between consecutive zeros, or equivalently, a normalized histogram of gaps between consecutive zeros. We work with the following cuspidal newforms:
\begin{center}
\begin{tabular}{c|c|c}
    LMFDB Label & Type & Group \\
    \hline
    \texttt{11.2.a.a}, \texttt{5.8.a.a}, \texttt{7.4.a.a} & $\chi_f$ principal, even twists & $\SpO(2N)$\\
    \texttt{11.2.a.a}, \texttt{5.8.a.a}, \texttt{7.4.a.a} & $\chi_f$ principal, odd twists & $\SpO(2N+1)$\\
    \texttt{3.7.b.a} & self-CM & $\USp(2N)$\\
    \texttt{13.2.e.a}, \texttt{11.7.b.b} & generic & $\U(N)$\\
\end{tabular}
\end{center}

In the following, we consider the `first' eigenvalues of random matrices $\SpO(\text{odd})$ and $\SpO(\text{even})$ whose characteristic polynomials are evaluated at or near 1. By `first,' we mean those eigenvalues closest to 1 on the unit circle. Note the eigenvalues of random matrices from $\SpO(\text{odd})$ are all going to be zero. 

\subsection{Families with orthogonal symmetry.}
\begin{center}
\begin{figure}[htpb]	
	\begin{tabular}{c c}
		Lowest zeros (even twists) & Lowest zeros (odd twists) \\
	\includegraphics[scale=0.45]{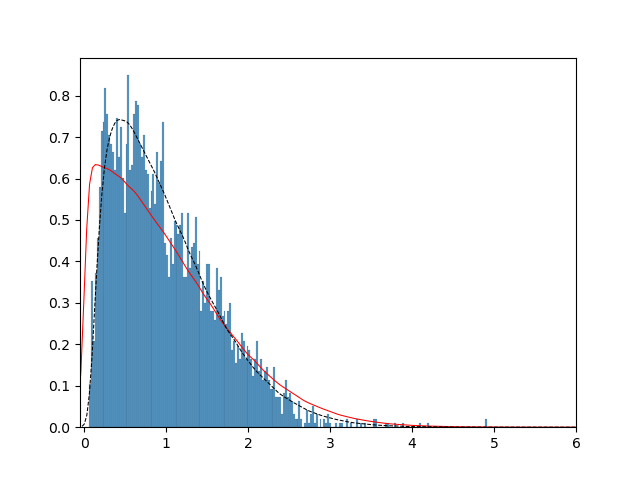} & \includegraphics[scale=0.45]{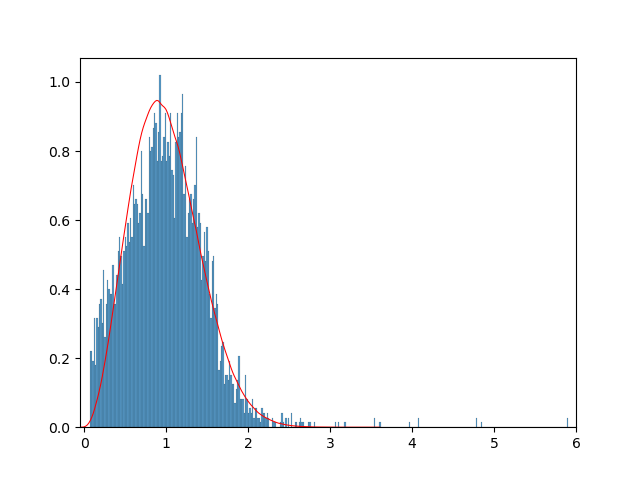}
	\end{tabular}
	\caption{The left histogram shows the distribution of lowest zeros for 4,500 even twists of \texttt{11.2.a.a}; the red curve (left) shows the distribution of the first eigenvalues of 1,000,000 randomly generated $\text{SO}(20)$ matrices with characteristic polynomial evaluated at 1; the black dotted curve (left) is the same distribution but with excision. We varied the excision threshold numerically to obtain the optimal fit. The right histogram shows the distribution of lowest zeros for 4,563 odd twists of \texttt{11.2.a.a}, and the red line (right) shows the distribution of the first eigenvalues of 1,000,000 randomly generated $\text{SO}(21)$ matrices with characteristic polynomial evaluated near 1. The data have been normalized to have mean 1.}
 \label{tab:zeros-mf11w2aa}
\end{figure}	
\end{center}

We numerically computed the lowest-lying zeros of $\mathcal{F}_f^{+}(X)$ with $f \in S_{2\cdot1}^{\text{new}}(11,\text{principal})$, which has sign $\epsilon_f=+1$. As expected, we obtain the repulsion from the origin for both even and odd twists and hence, the model requires the cut-off value. Recall we chose those twists with positive fundamental discriminants. If we choose twists with negative fundamental discriminants, we still recover repulsion from the origin. Our results combined with the results in \cite{DHKMS12} means that regardless of if we range over twists with negative or positive discriminants, we still recover repulsion from the origin. In addition, observe in Figure \ref{tab:zeros-mf11w2aa} the rather pronounced repulsion from the origin for lowest-lying zeros of odd twists of \texttt{11.2.a.a}. Since the distribution of the lowest-lying zeros agrees most with that of the excised matrices (dotted line) in Figure \ref{tab:zeros-mf11w2aa}, this verifies the necessity of creating such a model.
The lowest-lying zeros of odd twists do not vanish; this disagrees with eigenvalues of random matrices from $\SpO(\text{odd})$ which do, in fact, vanish. In a sense, the odd twists `force' a particular distribution which does not look as natural as that of even twists. In fact, we see an attraction toward the origin in the right histogram of Figure \ref{tab:zeros-mf11w2aa}. In the main paper \cite{BBJMNSY}, we present data for the form \texttt{5.4.a.a} which suggests this behaviour does not depend on the sign of the given form.

\begin{center}
\begin{figure}[htpb]
	\begin{tabular}{c c}
	Lowest zeros (even)
 & Non-vanishing lowest zeros (odd) \\
	\includegraphics[scale=0.45]{Images/Normalized_mf5w8aa_50000-sign+1_SO16.png}
 & \includegraphics[scale=0.45]{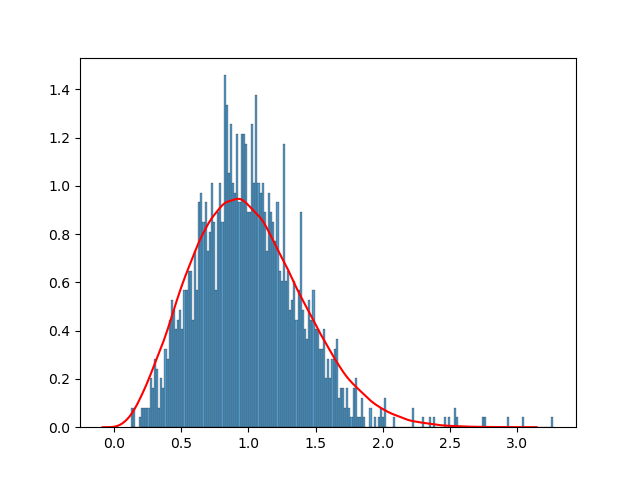}	
	\end{tabular}
	\caption{The left histogram shows lowest-lying zeros of 1,380 even twists of \texttt{5.8.a.a}; the red curve (left) shows the distribution of the first eigenvalues of 1,000,000 randomly generated $\SpO(16)$ matrices with characteristic polynomial evaluated at 1 without excision. The right histogram shows the non-vanishing lowest-lying zeros of 1,560 odd twists of \texttt{5.8.a.a}; the red curve (right) shows the distribution of the first eigenvalues of 1,000,000 randomly generated $\SpO(17)$ matrices with characteristic polynomial evaluated near 1. The data have been normalized to have mean 1.}
 \label{tab:table-mf5w8aa}
\end{figure}	
\end{center}

We numerically computed the lowest-lying zeros of $\mathcal{F}_f^{+}(X)$ with $f \in S_{2\cdot 4}^{\text{new}}(5,\text{principal})$, and the results are shown in Figure \ref{tab:table-mf5w8aa}. As predicted, the distribution of the lowest-lying zeros of even twists matches the eigenvalues of random matrices in $\SpO(\text{even})$ with characteristic polynomials evaluated. We also get agreement between the lowest-zeros of odd twists and the eigenvalues of random matrices in $\SpO(\text{odd})$ with characteristic polynomials evaluated at 1 as both vanish. As predicted, the distribution of the non-vanishing lowest-lying zeros of odd twists matches the eigenvalues of random matrices in $\SpO(\text{odd})$ with characteristic polynomials evaluated near 1, respectively. As opposed to the level 11, weight 2, principal form with even sign, the level 5, weight 8, principal form has odd sign $\epsilon_f=-1$ which indicates the odd twists do not `force' any $\SpO(\text{odd})$ behavior. Again, there is no observed repulsion from the origin and hence no need for a cut-off value. However, this phenomena might be caused by small data and so the repulsion is perhaps too small to be noticeable so quickly with the height increasing.

\begin{center}
\begin{figure}[htpb]	
\begin{tabular}{c c}
Lowest zeros (even twists) & Lowest zeros (odd twists) \\
\includegraphics[scale=0.323]{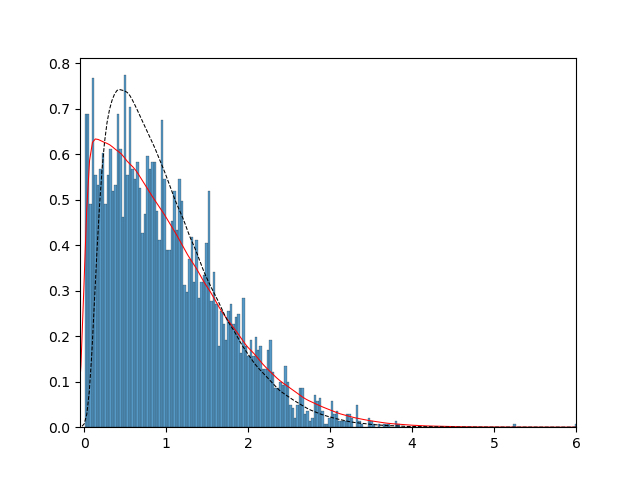} & \includegraphics[scale=0.45]{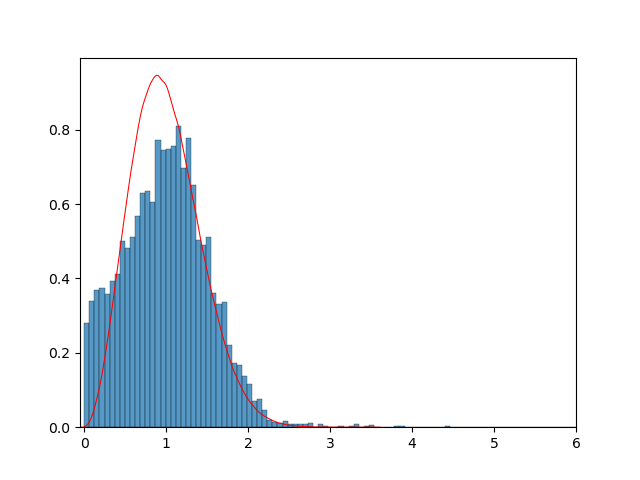}
\end{tabular}
\caption{The left histogram shows lowest-lying zeros of 4,684 even twists of \texttt{7.4.a.a}, and the red curve (left) shows the distribution of the first eigenvalues of 1,000,000 randomly generated $\SpO(20)$ matrices with characteristic polynomial evaluated at 1, and the black dotted curve is the same distribution but with excision with cutoff threshold the same as \texttt{11.2.a.a}. The right histogram shows non-vanishing lowest-lying zeros of 1,590 odd twists of \texttt{7.4.a.a}, and the red curve (right) shows the distribution of the first eigenvalues of 1,000,000 randomly generated $\SpO(21)$ matrices with characteristic polynomial evaluated near 1. The data have been normalized to have mean 1.}
\label{tab:zeros-mf7w4aa}
\end{figure}	
\end{center}

The sign of $f \in S_{2\cdot2}^{\text{new}}(7,\text{principal})$ is $\epsilon_f=+1$. As shown in the left histogram of Figure \ref{tab:zeros-mf7w4aa}, there is little to no discernible repulsion at the origin for the even twists. There is no need for excision for this form as shown by the disagreement between the excised distribution (dotted line) in Figure \ref{tab:zeros-mf7w4aa}; that is, the non-excised random matrix model describes the distribution well. This is expected given the heuristic proposed in \cite{DHKMS12}. However, on the right of Figure \ref{tab:zeros-mf7w4aa}, the red curve deviates from the data. In fact, we see the same attraction toward the origin for odd twists of \texttt{7.4.a.a} as that for odd twists of \texttt{11.2.a.a}. We notice both forms have even sign $\epsilon_f=+1$. This might indicate that attraction toward the origin for odd twists is dependent on the sign of the form being even. Rather than implement a cutoff value, one might develop a new model that accounts for this attraction by introducing a value that incorporates more first eigenvalues near the origin.

\subsection{Families with symplectic symmetry.} 
\begin{center}
\begin{figure}[htpb]	
	\begin{tabular}{c c c}
		Lowest zeros ($\Delta=+1$) & Lowest zeros ($\Delta=-1$)\\
	\includegraphics[scale=0.45]{Images/Normalized_mf3w7ba_50000_delta+1_USp20.png} & \includegraphics[scale=0.45]{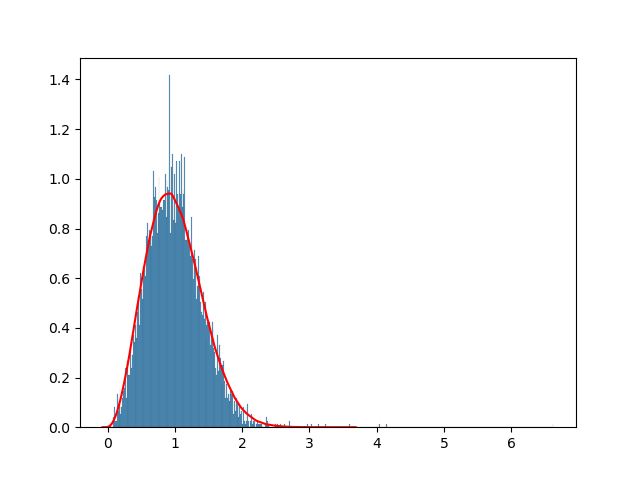} 
	\end{tabular}
	\caption{The left histogram shows the distribution of lowest-lying zeros for 5,450 twists of \texttt{3.7.b.a} with choice $\Delta =+1$, and the red line (left) shows the distribution of first eigenvalues of 1,000,000 randomly generated $\USp(20)$ matrices with characteristic polynomial evaluated at 1. The right histogram shows the distribution of lowest zeros for 5,720 twists of \texttt{3.7.b.a} with choice $\Delta =-1$, and the red line (right) shows the distribution of first eigenvalues of 1,000,000 randomly generated $\USp(20)$ matrices with characteristic polynomial evaluated at 1. The data have been normalized to have mean 1.}
 \label{tab:zeros-mf3w7ba}
\end{figure}	
\end{center}

As shown in Figure \ref{tab:zeros-mf3w7ba}, the lowest-lying zeros of $\mathcal{F}_f^{+}(X)$ for $f \in S_{7}^{\text{new}}(3,\text{self-CM})$ with sign $\epsilon_f=+1$ follow the predicted symplectic symmetry. The theory and the numerical results align as predicted.

\subsection{Families with unitary symmetry.}
\label{sec:families-with-unitary-symmetry}
\begin{center}
\begin{figure}[htpb]	
\begin{tabular}{c c}
Lowest zeros (twists) \\
\includegraphics[scale=0.55]{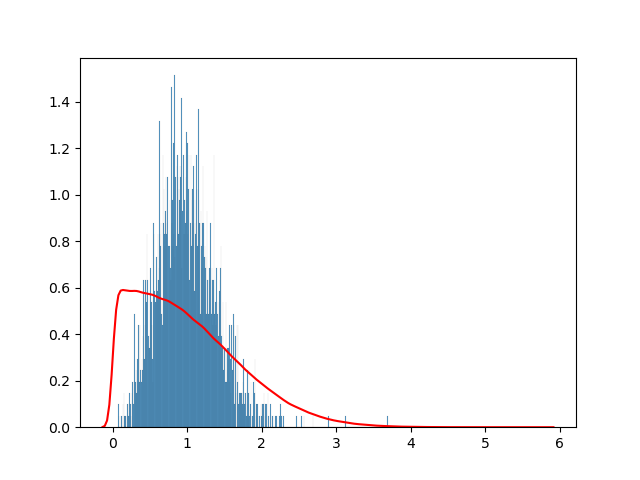} 
\end{tabular}
\caption{The histogram (blue) shows the distribution of the 2,860 lowest-lying zeros of the family of the generic form \texttt{11.7.b.b}; the red line shows the distribution of the eigenvalues of 1,000,000 random matrices from the unitary ensemble whose characteristic polynomials are evaluated at 1. The data have been normalized to have mean 1.}
\label{tab:zeros-mf11w7bb}
\end{figure}	
\end{center}

As shown in Figure \ref{tab:zeros-mf11w7bb}, the lowest-lying zeros of our family associated to the generic form \texttt{11.7.b.b} does not follow the predicted unitary distribution. The distribution of the low-lying zeros seems to match the distribution of the first eigenvalues of 1,000,000 numerically generated symplectic matrices. This means we recovered self-CM behavior from a generic form. Note the form \texttt{11.7.b.a} is, in fact, self-CM. Hence, a form with predicted unitary symmetry can have a different predicted symmetry under certain conditions—which are yet to be determined. The deviating behavior showcased in Figure \ref{tab:zeros-mf11w7bb} may be explained by the fact that the one-level density for the unitary ensemble showcases no oscillatory behavior as it equals 1. Hence, there is no possibility of extracting any arithmetic nuance. In particular, the data suggests certain generic forms that have unitary symmetry would restrict to have symmetry of another ensemble.

\begin{center}
\begin{figure}[htpb]	
\begin{tabular}{c c c}
Lowest zeros & Lowest zeros \\
\includegraphics[scale=0.33]{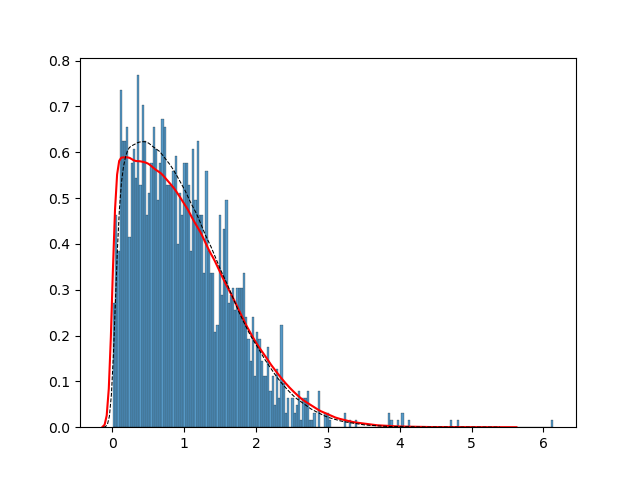} & \includegraphics[scale=0.45]{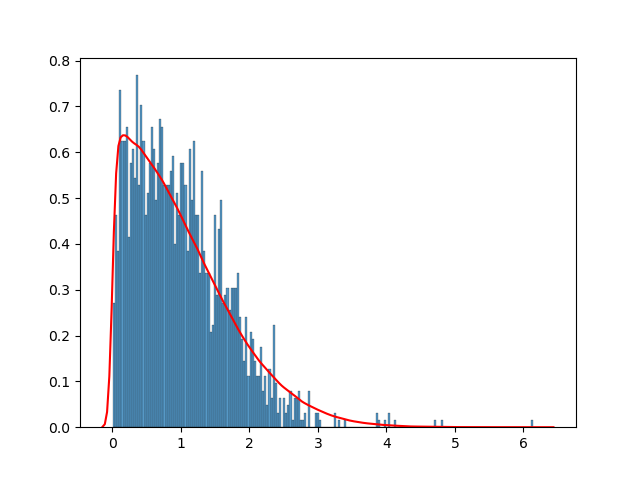} 
\end{tabular}
\caption{The left histogram shows the distribution of lowest-lying zeros for 2,097 twists of \texttt{13.2.e.a}; the red curve (left) shows the distribution of first eigenvalues of 1,000,000 randomly generated matrices from $\U(16)$ with characteristic polynomial evaluated at 1; the black dotted curve (left) is the same distribution but with excision. We varied the excision threshold numerically to obtain the optimal fit. The right histogram shows the same twists as the left histogram but has the red line (right) showing the distribution of first eigenvalues of 1,000,000 randomly generated matrices from $\SpO(16)$. The data have been normalized to have mean 1.}
 \label{tab:zeros-mf13w2ea}
\end{figure}	
\end{center}

The next generic form we considered was \texttt{13.2.e.a}, and the histograms are presented in Figure \ref{tab:zeros-mf13w2ea}. In \cite{BBJMNSY}, we took even and odd twists by setting $\psi_d(M)$ equal to either $+1$ or $-1$ to see if we recovered $\SpO(\text{even})$ and $\SpO(\text{odd})$ symmetry, respectively. We in fact did not recover $\SpO(\text{odd})$ symmetry, which is expected given that generic forms should not be influenced by parity of the sign. 
We transposed the first eigenvalues of 1,000,000 numerically generated special orthogonal (even) and unitary matrices, each normalized to have mean 1. This indistinguishability between the distributions presents difficulty for numerically determining the predicted ensemble for a generic form. We also see repulsion at the origin for \texttt{13.2.e.a}. The black dotted curve (left) is the distribution of random matrices from $\U(16)$ but with excision that was varied numerically to find the optimal fit. This indicates the excision present in the model should be extended to all weight 2 forms regardless of principality of the nebentype.

\newpage
\noindent


\end{document}